# Computing the Best Approximation Over the Intersection of a Polyhedral Set and the Doubly Nonnegative Cone

Ying Cui* Defeng Sun† Kim-Chuan Toh‡

**Abstract**

This paper introduces an efficient algorithm for computing the best approximation of a given matrix onto the intersection of linear equalities, inequalities and the doubly nonnegative cone (the cone of all positive semidefinite matrices whose elements are nonnegative). In contrast to directly applying the block coordinate descent (BCD) type methods, we propose an inexact accelerated (two-)block coordinate descent algorithm to tackle the four-block unconstrained nonsmooth dual program. The proposed algorithm hinges on the efficient semismooth Newton method to solve the subproblems, which have no closed form solutions since the original four blocks are merged into two larger blocks. The $O(1/k^2)$ iteration complexity of the proposed algorithm is established. Extensive numerical results over various large scale semidefinite programming instances from relaxations of combinatorial problems demonstrate the effectiveness of the proposed algorithm.

**Keywords.** semidefinite programming, doubly nonnegative cone, semismooth Newton, acceleration, complexity

**AMS subject classifications:** 90C06, 90C22, 90C25

## 1 Introduction

Let $\mathbb{X}$, $\mathbb{Y}$, and $\mathbb{Z}$ be three finite dimensional Euclidean spaces. Our aim in this paper is to solve the following convex optimization problem:

$$\begin{array}{ll} \underset{x\in\mathbb{X}}{\text{minimize}} & f(x)+\phi(x) \\ \text{subject to} & \mathcal{A}x=b, \quad g(x)\in\mathcal{C}, \quad x\in\mathcal{K}, \end{array} \tag{1}$$

where $f:\mathbb{X}\to\mathbb{R}$ is a smooth and strongly convex function, $\phi:\mathbb{X}\to(-\infty,+\infty]$ is a closed proper convex function, $\mathcal{A}:\mathbb{X}\to\mathbb{Y}$ is a linear operator, $b\in\mathbb{Y}$ is the given data, $\mathcal{C}\subseteq\mathbb{Z}$ and $\mathcal{K}\subseteq\mathbb{X}$ are two closed convex cones and $g:\mathbb{X}\to\mathbb{Z}$ is a smooth and $\mathcal{C}$-convex map for some closed convex set $\mathcal{C}$ (see, e.g., [34, Example 4']) satisfying

$$g\left(\alpha x_1+(1-\alpha)x_2\right)-\left[\,\alpha g(x_1)+(1-\alpha)g(x_2)\,\right]\in\mathcal{C}, \quad \forall\, x_1,x_2\in g^{-1}(\mathcal{C}),\ \forall\,\alpha\in(0,1).$$

*The Daniel J. Epstein Department of Industrial and Systems Engineering, University of Southern California, Los Angeles, California, USA. Email: yingcui@usc.edu

†Department of Applied Mathematics, The Hong Kong Polytechnic University, Hung Hom, Hong Kong. This research is partially supported by a start-up research grant from the Hong Kong Polytechnic University. Email: defeng.sun@polyu.edu.hk

‡Department of Mathematics and Institute of Operations Research and Analytics, National University of Singapore, Singapore. Email: mattohkc@nus.edu.sg

Our interest in solving (1) stemmed initially from a particular application pertaining to the best approximation problem over the intersection of affine spaces defined by linear equalities, inequalities and the doubly nonnegative cone, which takes the form of

$$\begin{array}{ll} \underset{X \in \mathbb{S}^n}{\text{minimize}} & \frac{1}{2} \| X - G \|^2 \\ \text{subject to} & \mathcal{A}X = b, \quad \mathcal{B}X \geq d, \quad X \in \mathbb{S}^n_+, \quad X \geq 0, \end{array} \tag{2}$$

where $\mathcal{A} : \mathbb{S}^n \to \mathbb{R}^{m_E}$, $\mathcal{B} : \mathbb{S}^n \to \mathbb{R}^{m_I}$ are linear operators, $G \in \mathbb{S}^n$, $b \in \mathbb{R}^{m_E}$, $d \in \mathbb{R}^{m_I}$ are given data, $\mathbb{S}^n$ and $\mathbb{S}^n_+$ are the cones of all $n \times n$ symmetric matrices and positive semidefinite matrices, respectively. Problem (1) also includes the following quadratically constrained quadratic program (which may include an additional nonsmooth regularization function $\phi$ in the objective) if $f$ is taken to be a strongly convex quadratic function and $g$ is set to be a map consisting of convex quadratic functions:

$$\begin{array}{ll} \underset{x \in \mathbb{X}}{\text{minimize}} & \frac{1}{2} \langle x, \mathcal{Q}_0 x \rangle + \langle c_0, x \rangle + \phi(x) \\ \text{subject to} & \frac{1}{2} \langle x, \mathcal{Q}_i x \rangle + \langle c_i, x \rangle + r_i \leq 0, \quad i = 1, \ldots, m, \quad x \in \mathcal{K}. \end{array}$$

Here $\mathcal{Q}_0 : \mathbb{X} \to \mathbb{X}$ is a self-adjoint positive definite linear operator, $\mathcal{Q}_i : \mathbb{X} \to \mathbb{X}$ are self-adjoint positive semidefinite linear operators, $c_0, c_i \in \mathbb{X}$ and $r_i \in \mathbb{R}$ are given data, all for $i = 1, \ldots, m$.

In this paper, we take a dual approach to solve (1) based on the observation that the dual program, with the form of

$$\underset{w}{\text{minimize}} \quad h(w) + \sum_{i=1}^{4} \varphi_i(w_i), \quad w \triangleq (w_1, w_2, w_3, w_4), \tag{3}$$

is an unconstrained convex nonsmooth problem, where $h$ is a convex differentiable function whose gradient is Lipschitz continuous, and $\varphi_1, \ldots, \varphi_4$ are proper closed convex functions; see section 2 for the detailed derivations of this dual formulation. In fact, besides being the dual program of (1), many optimization problems themselves have the form of (3), such as the robust principle component analysis [46] and the robust matrix completion problem [21].

Naturally, one may consider the block coordinate (gradient) descent (BCD) method, whose computational complexity is at best $O(1/k)$, to solve the four block unconstrained problem (3). See the papers [43, 42, 5, 33, 18] and a recent survey [47] for extensive discussions of this method. A key factor for determining the efficiency of the BCD method is the number of blocks that are updated sequentially during one iteration, since there is always a trade-off between such a number of blocks and the difficulty for solving the subproblems for each block. One may notice that solving the conic program (2) is different in nature from solving those problems arising in computational statistics and machine learning. The properties of the latter problems, such as the low computational cost of calculating one component of the gradient or solving one subproblem, and the need for only low accuracy solutions, are conducive for the efficient implementations of multi-block (usually at least hundreds-of-block) coordinate descent method. However, with the focus of solving the conic program (2) to a higher accuracy, we have the following issues to resolve:

- Treating (3) as a single block problem and solving all the variables simultaneously is difficult due to the degeneracy of this problem.

- Directly applying the four-block BCD method is inefficient as it potentially will need many more iteration cycles compared with those of three or fewer blocks. (This observation is indeed confirmed by the numerical experiments in section 5.) Within each iteration cycle, the algorithm may involve a computationally intensive step such as the projection onto the positive semidefinite cone.

To address the above issues, we propose to divide the four variables $(w_1, \ldots, w_4)$ into two groups and solve the problem (3) via a two-block inexact majorized accelerated BCD method. The subproblems with regard to each group may be non-degenerate and relatively easy to solve by Newton-type methods. With the presence of *four* separable nonsmooth functions, this algorithm is rather different from the inexact accelerated BCD method proposed in [39] to solve the doubly nonnegative best approximation problem. In [39], an additional regularization term related to the linear inequalities is added, which leads to only *two* nonsmooth separable functions in the corresponding dual program. Furthermore, by using Danskin's Theorem [12, Theorem 10.2.1], one of the nonsmooth blocks can be solved implicitly and the accelerated proximal gradient method initiated by Nesterov [27] is thus applicable to the resulting problem.

The key ingredient of our proposed algorithm is a combination of the inexactness, blockwise-updating, as well as the Nesterov-type acceleration technique. The attractive theoretical complexity of the acceleration technique for gradient-type methods and its good performance in practice has spurred much of the recent research to further investigate this technique. In particular, it has been extended to the BCD-type method with a random order in updating the blocks for smooth problems in [28] and nonsmooth problems in [23, 24, 1, 29]; and to a BCD method which is added with a large proximal term that is proportional to the number of blocks in [13, 14]. When the function $h$ in (3) is a least-squares quadratic function, the $O(1/k^2)$ complexity of a special majorized accelerated BCD algorithm is established in [7]. Though all of the above works have presented important theoretical progress, the numerical experiments in this paper indicate that their practical performance can significantly lag behind our proposed algorithm for solving the best approximation problem (2).

To summarize, the main contributions of our paper are as follows.

- We design an inexact two-block accelerated BCD method for solving (3) where $h$ is a general smooth coupled function that is not necessarily quadratic. A key feature of the method is the inexactness framework that allows us to apply the two-block based method to solve the 4-block problem (3) by viewing the 4 blocks as two larger blocks each containing 2 smaller blocks. Theoretically, we establish the $O(1/k^2)$ complexity of the proposed inexact accelerated BCD method. To achieve good practical performance, we suggest a proper way to merge several variables into one block so that only small proximal terms are added to the BCD subproblems. We also address the important issue of finding efficient algorithms to solve the BCD subproblems.

- Numerically, we provide an efficient solver based on the proposed method for solving the important class of best approximation problems of the form (2) that involves the positive semidefinite cone constraint and a large number of linear equality and inequality constraints. Our experiments on a large number of data instances from the Biq Mac Library maintained by A. Wiegele demonstrate that our solver is at least 3 to 4 times faster than other (accelerated) BCD-type methods.

The rest of the paper is organized as follows. In the next section, we derive the dual form of the problem (1) and propose the inexact majorized accelerated BCD method. Section 3 is devoted to the analysis of $O(1/k^2)$ iteration complexity of the proposed algorithm. In section 4, we describe the implementation of this inexact framework for solving the dual program (3). Newton-type algorithms for solving the subproblems are also discussed in this section. Numerical results are reported in section 5, where we show the effectiveness of our proposed algorithm via comparison with several variants of the BCD-type methods. We conclude our paper in section 6.

# 2 Formulation of the Dual Problem and the Algorithmic Framework

By introducing an auxiliary variable $\tilde{x} = x$ and reformulating (1) as

$$\begin{array}{ll} \underset{x,\tilde{x}\in\mathbb{X}}{\text{minimize}} & f(x) + \phi(\tilde{x}) \\ \text{subject to} & \mathcal{A}x = b, \quad g(x) \in \mathcal{C}, \quad x \in \mathcal{K}, \quad x = \tilde{x}, \end{array}$$

we derive the following Lagrangian function associated with the dual variable $(y, \lambda, s, z) \in \mathbb{Y} \times \mathbb{Z} \times \mathbb{X} \times \mathbb{X}$:

$$\mathcal{L}(x, \tilde{x}; y, \lambda, s, z) \triangleq f(x) + \phi(\tilde{x}) - \langle\, y,\, \mathcal{A}x - b\,\rangle - \langle\, \lambda,\, g(x)\,\rangle - \langle\, s,\, x\,\rangle - \langle\, z,\, x - \tilde{x}\,\rangle,$$

which leads to the dual program:

$$\begin{array}{ll} \underset{y,\lambda,s,z}{\text{maximize}} & \psi(\mathcal{A}^* y + s + z, \lambda) + \langle b, y\rangle - \phi^*(-z) \\ \text{subject to} & \lambda \in \mathcal{C}^*, \quad s \in \mathcal{K}^*. \end{array} \tag{4}$$

Here $\mathcal{A}^*$ is the adjoint of $\mathcal{A}$, $p^*$ is the conjugate function of $p$, $\mathcal{C}^*$ and $\mathcal{K}^*$ are the dual cones of $\mathcal{C}$ and $\mathcal{K}$, and the function $\psi : \mathbb{X} \times \mathcal{C}^* \to \mathbb{R}$ is defined as

$$\psi(w, \lambda) \triangleq \inf_{x\in\mathbb{X}} \{\, f(x) - \langle w, x\rangle - \langle\, \lambda,\, g(x)\,\rangle \,\}, \quad (w, \lambda) \in \mathbb{X} \times \mathcal{C}^*.$$

Since $g$ is assumed to be $\mathcal{C}$-convex, the term $-\langle\, \lambda,\, g(x)\,\rangle$ is convex with respect to $x$ for $\lambda \in \mathcal{C}^*$. The optimal solution of the above problem is thus a singleton by the strong convexity of $f$. In addition, the function $\psi$ is concave and continuously differentiable with Lipschitz continuous gradient [12, Theorem 10.2.1]. It therefore follows that the dual problem (4) is a special case of (3).

In what follows, we introduce an inexact majorized accelerated two-block coordinate descent method (imABCD) to solve (3). By grouping the 4 blocks of variables $(w_1, w_2, w_3, w_4)$ into two larger blocks, we can express (3) in the following 2-block format:

$$\underset{w,u,v}{\text{minimize}} \quad \theta(w) \triangleq h(w) + p(u) + q(v), \quad w \equiv (u, v),\ u \equiv (u_1, u_2) \in \mathbb{U},\ v \equiv (v_1, v_2) \in \mathbb{V}, \tag{5}$$

where $p(u) = \varphi_1(w_1) + \varphi_2(w_2)$ and $q(v) = \varphi_3(w_3) + \varphi_4(w_4)$ are proper closed convex functions. Since $\nabla h$ is assumed to be globally Lipschitz continuous, there exist two self-adjoint positive semidefinite linear operators $\mathcal{Q}$ and $\widehat{\mathcal{Q}} : \mathbb{U} \times \mathbb{V} \to \mathbb{U} \times \mathbb{V}$ such that

$$\left\{ \begin{array}{l} h(w) \,\geq\, h(w') + \langle \nabla h(w'), w - w'\rangle + \dfrac{1}{2}\, \|w - w'\|^2_{\mathcal{Q}}\,, \\[2ex] h(w) \,\leq\, \widehat{h}(w; w') \triangleq h(w') + \langle \nabla h(w'), w - w'\rangle \,+ \dfrac{1}{2}\, \|w - w'\|^2_{\widehat{\mathcal{Q}}}\,, \end{array} \right. \qquad \forall\, w, w' \in \mathbb{U} \times \mathbb{V}. \tag{6}$$

The operators $\mathcal{Q}$ and $\widehat{\mathcal{Q}}$ may be decomposed into the following $2 \times 2$ block structures as

$$\mathcal{Q}w \equiv \begin{pmatrix} \mathcal{Q}_{11} & \mathcal{Q}_{12} \\ \mathcal{Q}_{12}^* & \mathcal{Q}_{22} \end{pmatrix} \begin{pmatrix} u \\ v \end{pmatrix}, \quad \widehat{\mathcal{Q}}w \equiv \begin{pmatrix} \widehat{\mathcal{Q}}_{11} & \widehat{\mathcal{Q}}_{12} \\ \widehat{\mathcal{Q}}_{12}^* & \widehat{\mathcal{Q}}_{22} \end{pmatrix} \begin{pmatrix} u \\ v \end{pmatrix}, \quad w \in \mathbb{U} \times \mathbb{V},$$

where $\mathcal{Q}_{11}, \widehat{\mathcal{Q}}_{11} : \mathbb{U} \to \mathbb{U}$ and $\mathcal{Q}_{22}, \widehat{\mathcal{Q}}_{22} : \mathbb{V} \to \mathbb{V}$ are self-adjoint positive semidefinite linear operators, and $\mathcal{Q}_{12}, \widehat{\mathcal{Q}}_{12} : \mathbb{V} \to \mathbb{U}$ are two linear mappings whose adjoints are given by $\mathcal{Q}_{12}^*$ and $\widehat{\mathcal{Q}}_{12}^*$, respectively. The following assumption is made in the subsequent discussions.

**Assumption 1.** There exist two self-adjoint positive semidefinite linear operators $\mathcal{D}_1 : \mathbb{U} \to \mathbb{U}$ and $\mathcal{D}_2 : \mathbb{V} \to \mathbb{V}$ such that

$$\widehat{\mathcal{Q}} \;=\; \mathcal{Q} + \operatorname{Diag}(\mathcal{D}_1, \mathcal{D}_2).$$

Furthermore, $\widehat{\mathcal{Q}}$ satisfies that $\widehat{\mathcal{Q}}_{11} = \mathcal{Q}_{11} + \mathcal{D}_1 \succ 0$ and $\widehat{\mathcal{Q}}_{22} = \mathcal{Q}_{22} + \mathcal{D}_2 \succ 0$.

Below is our proposed algorithm for solving problem (3) via the two-block program (5).

---

**imABCD:** An **i**nexact **m**ajorized **A**ccelerated **B**lock **C**oordinate **D**escent algorithm for (5)

---

**Initialization**. Choose an initial point $(u^1, v^1) = (\widetilde{u}^0, \widetilde{v}^0) \in \operatorname{dom} p \times \operatorname{dom} q$ and a nonnegative non-increasing sequence $\{\varepsilon_k\}$. Let $t_1 = 1$. Perform the following steps in each iteration for $k \geq 1$.

**Step 1.** Compute

$$\begin{cases} \widetilde{u}^k = \operatorname*{argmin}_{u \in \mathbb{U}} \Big\{ p(u) + \widehat{h}(u, v^k; w^k) + \langle \delta_u^k, u \rangle \Big\}, \\ \widetilde{v}^k = \operatorname*{argmin}_{v \in \mathbb{V}} \Big\{ q(v) + \widehat{h}(\widetilde{u}^k, v; w^k) + \langle \delta_v^k, v \rangle \Big\} \end{cases} \tag{7}$$

such that $(\delta_u^k, \delta_v^k) \in \mathbb{U} \times \mathbb{V}$ satisfies $\max\Big\{ \|\widehat{\mathcal{Q}}_{11}^{-1/2} \delta_u^k\|, \, \|\widehat{\mathcal{Q}}_{22}^{-1/2} \delta_v^k\| \Big\} \leq \varepsilon_k$. Denote $\widetilde{w}^k = (\widetilde{u}^k, \widetilde{v}^k)$.

**Step 2.** Compute
$$\begin{cases} t_{k+1} = \dfrac{1}{2}\left(1 + \sqrt{1 + 4t_k^2}\right), \\ w^{k+1} = \widetilde{w}^k + \dfrac{t_k - 1}{t_{k+1}}\left(\widetilde{w}^k - \widetilde{w}^{k-1}\right). \end{cases}$$

---

The above imABCD algorithm can be taken as an accelerated as well as an inexact version of the alternating minimization method. When $\varepsilon_k \equiv 0$ for all $k \geq 0$, the proposed algorithm reduces to an exact version of the majorized accelerated block coordinate descent (mABCD) method. Since Nesterov's acceleration technique is able to improve the complexity of the gradient-type method from $O(1/k)$ to $O(1/k^2)$, it is interesting for us to investigate whether this acceleration technique can be extended to the two-block coordinate descent method without random selection of the updating blocks. In fact, extensive numerical experiments in the existing literature indicate that the acceleration technique may substantially improve the efficiency of the BCD algorithm, see, e.g., the numerical comparison in [39]. The study on this subject is thus critical for understanding the reasons behind this phenomenon.

## 3 The $O(1/k^2)$ Complexity of the Objective Values

In order to simplify the subsequent discussions, we introduce the positive semidefinite operator $\mathcal{H} : \mathbb{U} \times \mathbb{V} \to \mathbb{U} \times \mathbb{V}$ defined by

$$\mathcal{H} \triangleq \mathrm{Diag}\,(\,\mathcal{D}_1\,,\,\mathcal{D}_2 + \mathcal{Q}_{22}\,). \tag{8}$$

We also write $\Omega$ as the optimal solution set of (3). We start by presenting the the following simple lemma concerning the properties of the sequence $\{t_k\}$.

**Lemma 1.** The sequence $\{t_k\}_{k\geq 1}$ generated by the imABCD algorithm satisfies the following properties:

(a) $1 - \dfrac{1}{t_{k+1}} = \dfrac{t_k^2}{t_{k+1}^2}$ .   (b) $\dfrac{k+1}{2} \leq t_k \leq \dfrac{5}{8}k + \dfrac{3}{8} \leq k$ .

(c) Given any $w^+, w'$ in $\mathbb{U} \times \mathbb{V}$. Consider $w = (1 - \frac{1}{t_k})w^+ + \frac{1}{t_k}w'$. Then

$$t_k^2\,[\,\theta(w) - \theta(w')\,] \;\leq\; t_{k-1}^2\,[\,\theta(w^+) - \theta(w')\,]. \tag{9}$$

*Proof.* By noting that $t_{k+1}^2 - t_{k+1} = t_k^2$, the property (a) can be obtained directly. The property (b) can be derived from the inequalities

$$t_{k+1} = \frac{1 + 2t_k\sqrt{1 + 1/(4t_k^2)}}{2} \leq \frac{1 + 2t_k(1 + 1/(8t_k^2))}{2} = \frac{1}{2} + t_k + \frac{1}{8t_k} \leq \frac{5}{8} + t_k \leq \frac{5k}{8} + t_1 = \frac{5k}{8} + 1$$

and

$$t_{k+1} = \frac{1 + \sqrt{1 + 4t_k^2}}{2} \geq \frac{1 + 2t_k}{2} \geq \frac{k + 2t_1}{2} = \frac{k+2}{2}\,.$$

(c) From the convexity of $\theta$, we have that

$$t_k^2\,\theta(w) \;\leq\; (t_k^2 - t_k)\,\theta(w^+) + t_k\,\theta(w') \;=\; t_{k-1}^2\,\theta(w^+) + (t_k^2 - t_{k-1}^2)\,\theta(w').$$

From here, we get the desired inequality. □

In the following, we shall first provide the $O(1/k^2)$ complexity of the ABCD method with subproblems being solved exactly, and then extend the analysis to the inexact case.

### 3.1 The case where the subproblems are being solved exactly

The analysis in this subsection is partially motivated by the recent paper [7] in which the authors consider the $O(1/k^2)$ complexity of an accelerated BCD method for (3) where $h$ is a special least-squares quadratic function. Here we extend this nice result to a more general setting where $h$ is only required to be a smooth function.

The lemma below shows an important property of the objective values for the sequence generated by the mABCD algorithm, which is essential to prove the main global complexity result.

**Lemma 2.** Suppose that Assumption 1 holds. Let the sequences $\{\widetilde{w}^k\} \triangleq \{(\widetilde{u}^k, \widetilde{v}^k)\}$ and $\{w^k\} = \{(u^k, v^k)\}$ be generated by the mABCD algorithm. Then for any $k \geq 1$, it holds that

$$\theta(\widetilde{w}^k) - \theta(w) \;\leq\; \frac{1}{2}\,\|\, w - w^k \,\|_{\mathcal{H}}^2 - \frac{1}{2}\,\|\, w - \widetilde{w}^k \,\|_{\mathcal{H}}^2, \quad \forall\, w \in \mathbb{U} \times \mathbb{V}.$$

In particular, if $w^*$ is an optimal solution of (3), then

$$\|\, w^* - \widetilde{w}^k \|_{\mathcal{H}} \quad \leq \; \|\, w^* - w^k \,\|_{\mathcal{H}}.$$

*Proof.* By applying the optimality condition to the subproblems in (7), we derive that

$$\begin{cases} 0 \in \; \partial p(\widetilde{u}^k) + \nabla_u h(w^k) + \widehat{\mathcal{Q}}_{11}(\widetilde{u}^k - u^k)\;, \\ 0 \in \; \partial q(\widetilde{v}^k) + \nabla_v h(w^k) + \mathcal{Q}_{12}^*(\widetilde{u}^k - u^k) + \widehat{\mathcal{Q}}_{22}(\widetilde{v}^k - v^k)\;. \end{cases}$$

Thus, it follows from the convexity of $p(\cdot)$ and $q(\cdot)$ that

$$\begin{cases} p(u) \;\geq\; p(\widetilde{u}^k) + \Big\langle u - \widetilde{u}^k,\; -\nabla_u h(w^k) - \widehat{\mathcal{Q}}_{11}(\widetilde{u}^k - u^k)\Big\rangle\,, \quad \forall\, u \in \mathbb{U}, \\ q(v) \;\geq\; q(\widetilde{v}^k) + \Big\langle v - \widetilde{v}^k,\; -\nabla_v h(w^k) - \mathcal{Q}_{12}^*(\widetilde{u}^k - u^k) - \widehat{\mathcal{Q}}_{22}(\widetilde{v}^k - v^k)\Big\rangle\,, \quad \forall\, v \in \mathbb{V}. \end{cases} \tag{10}$$

Based on the inequalities in (6) that

$$\begin{cases} h(\widetilde{w}^k) \;\leq\; h(w^k) + \big\langle \nabla h(w^k), \widetilde{w}^k - w^k \big\rangle + \dfrac{1}{2}\,\|\, \widetilde{w}^k - w^k \,\|_{\widehat{\mathcal{Q}}}^2\;, \\ h(w) \;\geq\; h(w^k) + \big\langle \nabla h(w^k), w - w^k \big\rangle + \dfrac{1}{2}\,\|\, w - w^k \,\|_{\mathcal{Q}}^2\;, \end{cases}$$

we get

$$h(w) - h(\widetilde{w}^k) \;\geq\; \Big\langle \nabla h(w^k), w - \widetilde{w}^k \Big\rangle + \frac{1}{2}\,\|\, w - w^k \,\|_{\mathcal{Q}}^2 - \frac{1}{2}\,\|\, \widetilde{w}^k - w^k \,\|_{\widehat{\mathcal{Q}}}^2\;. \tag{11}$$

By the Cauchy-Schwarz inequality, we also have

$$\begin{aligned} &2\,\big\langle\, \widetilde{u}^k - u,\; \mathcal{Q}_{12}\,(\widetilde{v}^k - v^k)\,\big\rangle \;=\; 2\,\Big\langle \mathcal{Q}\,(\widetilde{w}^k - w),\; \begin{pmatrix} 0 \\ \widetilde{v}^k - v^k \end{pmatrix} \Big\rangle - 2\,\langle\, \mathcal{Q}_{22}\,(\widetilde{v}^k - v), \widetilde{v}^k - v^k \,\rangle \\ &\leq\; \big(\,\|\, \widetilde{w}^k - w \,\|_{\mathcal{Q}}^2 + \|\, \widetilde{v}^k - v^k \,\|_{\mathcal{Q}_{22}}^2\,\big) - \big(\,\|\, \widetilde{v}^k - v \,\|_{\mathcal{Q}_{22}}^2 + \|\, \widetilde{v}^k - v^k \,\|_{\mathcal{Q}_{22}}^2 - \|\, v^k - v \,\|_{\mathcal{Q}_{22}}^2\,\big) \\ &=\; \|\, \widetilde{w}^k - w \,\|_{\mathcal{Q}}^2 + \|\, v^k - v \,\|_{\mathcal{Q}_{22}}^2 - \|\, \widetilde{v}^k - v \,\|_{\mathcal{Q}_{22}}^2. \end{aligned} \tag{12}$$

Summing up the inequalities (10) and (11) and substituting the resulting inequality into (12), we obtain

$$\begin{aligned} &\quad 2\,\big(\,\theta(w) - \theta(\widetilde{w}^k)\,\big) \\ &\geq\; \Big(\,\|\, w - w^k \,\|_{\mathcal{Q}}^2 - \|\, \widetilde{w}^k - w^k \,\|_{\widehat{\mathcal{Q}}}^2\Big) - 2\,\Big\langle w - \widetilde{w}^k, \widehat{\mathcal{Q}}\,(\widetilde{w}^k - w^k)\Big\rangle - 2\,\big\langle\, \widetilde{u}^k - u, \mathcal{Q}_{12}\,(\widetilde{v}^k - v^k)\,\big\rangle \\ &\geq\; \|\, w - w^k \,\|_{\mathcal{Q}}^2 - \|\, \widetilde{w}^k - w^k \,\|_{\widehat{\mathcal{Q}}}^2 - \Big(\|\, w - w^k \,\|_{\widehat{\mathcal{Q}}}^2 - \|\, w - \widetilde{w}^k \,\|_{\widehat{\mathcal{Q}}}^2 - \|\, \widetilde{w}^k - w^k \,\|_{\widehat{\mathcal{Q}}}^2\Big) \\ &\quad - \big(\,\|\, \widetilde{w}^k - w \,\|_{\mathcal{Q}}^2 + \|\, v^k - v \,\|_{\mathcal{Q}_{22}}^2 - \|\, \widetilde{v}^k - v \,\|_{\mathcal{Q}_{22}}^2\,\big) \;=\; \|\, w - \widetilde{w}^k \,\|_{\mathcal{H}}^2 - \|\, w - w^k \,\|_{\mathcal{H}}^2\;, \end{aligned}$$

where the last equality is due to Assumption 1. The stated inequality therefore follows. □

Based on the above lemma, we next show the $O(1/k^2)$ complexity for the sequence of objective values obtained by the mABCD algorithm.

**Theorem 2.** Suppose that Assumption 1 holds and the solution set $\Omega$ of the problem (3) is non-empty. Let $w^* \triangleq (u^*, v^*) \in \Omega$. Then the sequence $\{\widetilde{w}^k\} \triangleq \{(\widetilde{u}^k, \widetilde{v}^k)\}$ generated by the mABCD algorithm satisfies that

$$\theta(\widetilde{w}^k) - \theta(w^*) \le \frac{2\,\|\,\widetilde{w}^0 - w^*\,\|_{\mathcal{H}}^2}{(k+1)^2}\,, \quad \forall\, k \ge 1.$$

*Proof.* Letting $w = (1 - \frac{1}{t_k})\widetilde{w}^{k-1} + \frac{1}{t_k} w^*$ in Lemma 2, we derive that, for any $k \ge 2$,

$$t_k^2\theta(w) \,-\, t_k^2\theta(\widetilde{w}^k) \ge \frac{1}{2}\left\|(t_k - 1)\widetilde{w}^{k-1} + w^* - t_k\widetilde{w}^k\right\|_{\mathcal{H}}^2 - \frac{1}{2}\left\|(t_k - 1)\widetilde{w}^{k-1} + w^* - t_k w^k\right\|_{\mathcal{H}}^2.$$

By applying Lemma 1 (c) with $w^+ = \widetilde{w}^{k-1}$ and $w' = w^*$, we get

$$t_k^2\,[\,\theta(w) - \theta(w^*)\,] \le t_{k-1}^2[\,\theta(\widetilde{w}^{k-1}) - \theta(w^*)\,].$$

By combining the above inequalities and noting that $w^k = \widetilde{w}^{k-1} + \frac{t_{k-1}-1}{t_k}(\widetilde{w}^{k-1} - \widetilde{w}^{k-2})$, we get for $k \ge 2$,

$$\begin{aligned}
& t_k^2\,[\,\theta(\widetilde{w}^k) - \theta(w^*)\,] + \frac{1}{2}\,\|\,t_k\widetilde{w}^k - w^* - (t_k - 1)\widetilde{w}^{k-1}\,\|_{\mathcal{H}}^2 \\
\le\ & t_{k-1}^2\,[\,\theta(\widetilde{w}^{k-1}) - \theta(w^*)\,] + \frac{1}{2}\,\|\,t_{k-1}\widetilde{w}^{k-1} - w^* - (t_{k-1} - 1)\widetilde{w}^{k-2}\,\|_{\mathcal{H}}^2.
\end{aligned}$$

By applying Lemma 2 again for $k = 1$ and $w = w^*$, we get

$$\theta(\widetilde{w}^1) - \theta(w^*) \le \frac{1}{2}\,\|w^1 - w^*\,\|_{\mathcal{H}}^2 - \frac{1}{2}\,\|\,\widetilde{w}^1 - w^*\,\|_{\mathcal{H}}^2 = \frac{1}{2}\,\|\,\widetilde{w}^0 - w^*\,\|_{\mathcal{H}}^2 - \frac{1}{2}\,\|\,\widetilde{w}^1 - w^*\,\|_{\mathcal{H}}^2\,.$$

It therefore follows that, for all $k \ge 1$,

$$\begin{aligned}
& t_k^2\,[\,\theta(\widetilde{w}^k) - \theta(w^*)\,] \,+\, \frac{1}{2}\,\|\,t_k\widetilde{w}^k - w^* - (t_k - 1)\widetilde{w}^{k-1}\,\|_{\mathcal{H}}^2 \\
\le\ & t_{k-1}^2\,[\,\theta(\widetilde{w}^{k-1}) - \theta(w^*)\,] \,+\, \frac{1}{2}\|\,t_{k-1}\widetilde{w}^{k-1} - w^* - (t_{k-1} - 1)\widetilde{w}^{k-2}\,\|_{\mathcal{H}}^2 \\
\le\ & \cdots \\
\le\ & t_1^2\,[\,\theta(\widetilde{w}^1) - \theta(w^*)\,] \,+\, \frac{1}{2}\,\|\,t_1\widetilde{w}^1 - w^* - (t_1 - 1)\widetilde{w}^0\,\|_{\mathcal{H}}^2 \le \frac{1}{2}\|\,\widetilde{w}^0 - w^*\,\|_{\mathcal{H}}^2\,.
\end{aligned}$$

The desired inequality of this theorem can thus be established since $t_k \ge \dfrac{k+1}{2}$ by Lemma 1(b). $\square$

### 3.2 The case where the subproblems are being solved inexactly

Theorem 2 shows the $O(1/k^2)$ complexity of the objective values for the two-block majorized accelerated BCD algorithm for solving (3). However, there seems to have no easy extension of its proof to problems with three or more blocks. In this section, we consider to allow for inexactness when solving the subproblems. The introduction of the inexactness is crucial for efficiently solving the multi-block problems. We note that such an idea has already been incorporated into different variants of the BCD and APG algorithms, see, e.g., [37, 20, 45, 41]. But the analyses therein are not applicable to our setting.

We characterize the decreasing property of the objective values for the imABCD algorithm in the lemma below. Its proof can be derived similarly as the proof of the exact case in Proposition 2 by applying the optimality conditions at the iteration point $(\widetilde{u}^k, \widetilde{v}^k)$. We omit the details here for brevity.

**Lemma 3.** Suppose that Assumption 1 holds. Let the sequences $\{\widetilde{w}^k\} \triangleq \{(\widetilde{u}^k, \widetilde{v}^k)\}$ and $\{w^k\} \triangleq \{(u^k, v^k)\}$ be generated by the imABCD algorithm. Then for any $k \geq 1$,

$$\theta(\widetilde{w}^k) - \theta(w) \leq \frac{1}{2} \| w - w^k \|_{\mathcal{H}}^2 - \frac{1}{2} \| w - \widetilde{w}^k \|_{\mathcal{H}}^2 + \varepsilon_k \| w - \widetilde{w}^k \|_{\mathrm{Diag}(\widehat{\mathcal{Q}}_{11}, \widehat{\mathcal{Q}}_{22})}, \quad \forall\, w \in \mathbb{U} \times \mathbb{V}.$$

For $k \geq 1$, we denote the exact solutions at the $(k+1)$-th iteration as

$$\overline{u}^k \triangleq \operatorname*{argmin}_{u \in \mathbb{U}} \left\{ p(u) + \widehat{h}(u, v^k; w^k) \right\}, \quad \overline{v}^k \triangleq \operatorname*{argmin}_{v \in \mathbb{V}} \left\{ q(v) + \widehat{h}(\overline{u}^k, v; w^k) \right\}. \tag{13}$$

For consistency, we set $(\overline{u}^0, \overline{v}^0) = (\widetilde{u}^0, \widetilde{v}^0) = (u^1, v^1)$. Since $\widehat{\mathcal{Q}}_{11}$ and $\widehat{\mathcal{Q}}_{22}$ are assumed to be positive definite, the above two problems admit unique solutions and thus, $\overline{u}^k$ and $\overline{v}^k$ are well defined for $k \geq 0$. The following lemma shows the gap between $(\overline{u}^k, \overline{v}^k)$ and $(\widetilde{u}^k, \widetilde{v}^k)$.

**Lemma 4.** Let the sequences $\{(\widetilde{u}^k, \widetilde{v}^k)\}$ and $\{(u^k, v^k)\}$ be generated by the imABCD algorithm, and $\{(\overline{u}^k, \overline{v}^k)\}$ be given by (13). For any $k \geq 1$, the following inequalities hold:

(a) $\| \overline{u}^k - u^* \|_{\widehat{\mathcal{Q}}_{11}}^2 \leq \| u^k - u^* \|_{\mathcal{D}_1}^2 + \| v^k - v^* \|_{\widehat{\mathcal{Q}}_{22}}^2 = \| w^k - w^* \|_{\mathcal{H}}^2$;

(b) $\| \widehat{\mathcal{Q}}_{11}^{1/2} (\widetilde{u}^k - \overline{u}^k) \| \leq \varepsilon_k, \quad \| \widehat{\mathcal{Q}}_{22}^{1/2} (\widetilde{v}^k - \overline{v}^k) \| \leq (1 + \sqrt{2})\, \varepsilon_k$;

(c) $\| \widetilde{w}^k - \overline{w}^k \|_{\mathcal{H}} \leq \sqrt{7}\, \varepsilon_k$.

*Proof.* (a) By applying the optimality conditions to the problems in (13), we deduce that

$$0 \in \partial p(\overline{u}^k) + \nabla_u h(w^k) + \widehat{\mathcal{Q}}_{11}(\overline{u}^k - u^k) \quad \text{and} \quad 0 \in \partial p(u^*) + \nabla_u h(w^*).$$

By the monotonicity of $\partial p$, we get

$$\left\langle \overline{u}^k - u^*, \, \nabla_u h(w^k) - \nabla_u h(w^*) + \widehat{\mathcal{Q}}_{11}(\overline{u}^k - u^k) \right\rangle \leq 0. \tag{14}$$

Since $\nabla h$ is assumed to be globally Lipschitz continuous, it is known from Clarke's Mean-Value Theorem [8, Proposition 2.6.5] that there exists a self-adjoint and positive semidefinite operator $W^k \in \operatorname{conv}\left\{ \partial^2 h(\, [w^{k-1}, w^k]\,) \right\}$ such that

$$\nabla h(w^k) - \nabla h(w^{k-1}) = W^k (w^k - w^{k-1}),$$

where the set $\operatorname{conv}\left\{ \partial^2 h(\, [w^{k-1}, w^k]\,) \right\}$ denotes the convex hull of all points in $\partial^2 h(z)$ for any $z \in [w^{k-1}, w^k]$. To proceed, we write $W^k \triangleq \begin{pmatrix} W_{11}^k & W_{12}^k \\ (W_{12}^k)^* & W_{22}^k \end{pmatrix}$, where $W_{11}^k : \mathbb{U} \to \mathbb{U}$, $W_{22}^k : \mathbb{V} \to \mathbb{V}$ are self-adjoint positive semidefinite linear operators and $W_{12}^k : \mathbb{U} \to \mathbb{V}$ is a linear operator. Since

$$\left\langle \begin{pmatrix} \bar{u}^k - u^* \\ v^k - v^* \end{pmatrix}, W^k \begin{pmatrix} \bar{u}^k - u^* \\ v^k - v^* \end{pmatrix} \right\rangle \geq 0 \quad \text{and} \quad \mathcal{Q} \preceq W^k \preceq \widehat{\mathcal{Q}},$$

we can derive that

$$
\begin{aligned}
&2\left\langle \overline{u}^k-u^*,\, \nabla_u h(w^k)-\nabla_u h(w^*) \right\rangle \;=\; 2\left\langle \overline{u}^k-u^*,\, W_{11}^k(u^k-u^*)+W_{12}^k(v^k-v^*) \right\rangle \\
&\quad\geq\; \Big(\|\,\overline{u}^k-u^*\,\|_{W_{11}^k}^2+\|\,u^k-u^*\,\|_{W_{11}^k}^2-\|\,\overline{u}^k-u^k\,\|_{W_{11}^k}^2\Big)-\Big(\|\,\overline{u}^k-u^*\,\|_{W_{11}^k}^2+\|\,v^k-v^*\,\|_{W_{22}^k}^2\Big) \\
&\quad\geq\; \|\,u^k-u^*\,\|_{\mathcal{Q}_{11}}^2-\|\,\overline{u}^k-u^k\,\|_{\widehat{\mathcal{Q}}_{11}}^2-\|\,v^k-v^*\,\|_{\widehat{\mathcal{Q}}_{22}}^2.
\end{aligned}
\tag{15}
$$

From (14) and (15), we may obtain that

$$
\begin{aligned}
0 \;\geq\;& \|\,u^k-u^*\,\|_{\mathcal{Q}_{11}}^2-\|\,\overline{u}^k-u^k\,\|_{\widehat{\mathcal{Q}}_{11}}^2-\|\,v^k-v^*\,\|_{\widehat{\mathcal{Q}}_{22}}^2+2\left\langle \overline{u}^k-u^*,\, \widehat{\mathcal{Q}}_{11}\,(\overline{u}^k-u^k)\right\rangle \\
=\;& \|\,u^k-u^*\,\|_{\mathcal{Q}_{11}}^2-\|\,\overline{u}^k-u^k\,\|_{\widehat{\mathcal{Q}}_{11}}^2-\|\,v^k-v^*\,\|_{\widehat{\mathcal{Q}}_{22}}^2+\|\,\overline{u}^k-u^*\,\|_{\widehat{\mathcal{Q}}_{11}}^2+\|\,\overline{u}^k-u^k\,\|_{\widehat{\mathcal{Q}}_{11}}^2-\|\,u^k-u^*\,\|_{\widehat{\mathcal{Q}}_{11}}^2 \\
=\;& \|\,\overline{u}^k-u^*\,\|_{\widehat{\mathcal{Q}}_{11}}^2-\|\,u^k-u^*\,\|_{\mathcal{D}_1}^2-\|\,v^k-v^*\,\|_{\widehat{\mathcal{Q}}_{22}}^2,
\end{aligned}
$$

which yields $\|\,\overline{u}^k-u^*\,\|_{\widehat{\mathcal{Q}}_{11}}^2 \;\leq\; \|\,u^k-u^*\,\|_{\mathcal{D}_1}^2+\|\,v^k-v^*\|_{\widehat{\mathcal{Q}}_{22}}^2$. This completes the proof of part (a).

(b) In order to obtain bounds for $\|\,\widehat{\mathcal{Q}}_{11}^{1/2}\,(\widetilde{u}^k-\overline{u}^k)\|$ and $\|\,\widehat{\mathcal{Q}}_{22}^{1/2}\,(\widetilde{v}^k-\overline{v}^k)\|$, we apply the optimality conditions to the problems in (13) at $(\overline{u}^k,\overline{v}^k)$ and $(\widetilde{u}^k,\widetilde{v}^k)$ to deduce that

$$
\left\langle \widehat{\mathcal{Q}}_{11}(\widetilde{u}^k-\overline{u}^k)+\delta_u^k,\, \widetilde{u}^k-\overline{u}^k\right\rangle \;\leq\; 0 \quad\text{and}\quad \left\langle \mathcal{Q}_{12}^*\,(\widetilde{u}^k-\overline{u}^k)+\widehat{\mathcal{Q}}_{22}\,(\widetilde{v}^k-\overline{v}^k)+\delta_v^k,\, \widetilde{v}^k-\overline{v}^k\right\rangle \;\leq\; 0.
$$

The first inequality implies that

$$
\|\,\widehat{\mathcal{Q}}_{11}^{1/2}\,(\widetilde{u}^k-\overline{u}^k)\,\| \;\leq\; \|\,\widehat{\mathcal{Q}}_{11}^{-1/2}\,\delta_u^k\,\| \;\leq\; \varepsilon_k.
$$

The second inequality yields that

$$
\begin{aligned}
\|\,\widetilde{v}^k-\overline{v}^k\,\|_{\widehat{\mathcal{Q}}_{22}}^2 \;&\leq\; \|\,\widehat{\mathcal{Q}}_{22}^{-1/2}\,\delta_v^k\,\|\,\|\,\widehat{\mathcal{Q}}_{22}^{1/2}\,(\widetilde{v}^k-\overline{v}^k)\,\|-\left\langle \mathcal{Q}_{12}^*\,(\widetilde{u}^k-\overline{u}^k),\widetilde{v}^k-\overline{v}^k\right\rangle \\
&\leq\; \|\,\widehat{\mathcal{Q}}_{22}^{-1/2}\,\delta_v^k\|\,\|\,\widehat{\mathcal{Q}}_{22}^{1/2}\,(\widetilde{v}^k-\overline{v}^k)\|+\frac{1}{2}\,\Big(\|\,\widetilde{u}^k-\overline{u}^k\,\|_{\widehat{\mathcal{Q}}_{11}}^2+\|\,\widetilde{v}^k-\overline{v}^k\,\|_{\widehat{\mathcal{Q}}_{22}}^2\Big).
\end{aligned}
$$

Hence,

$$
\|\,\widetilde{v}^k-\overline{v}^k\,\|_{\widehat{\mathcal{Q}}_{22}}^2 \;\leq\; 2\,\|\,\widehat{\mathcal{Q}}_{22}^{-1/2}\delta_v^k\,\|\,\|\,\widehat{\mathcal{Q}}_{22}^{1/2}(\widetilde{v}^k-\overline{v}^k)\,\|+\|\,\widetilde{u}^k-\overline{u}^k\,\|_{\widehat{\mathcal{Q}}_{11}}^2 \;\leq\; 2\varepsilon_k\,\|\,\widehat{\mathcal{Q}}_{22}^{1/2}(\widetilde{v}^k-\overline{v}^k)\,\|+\varepsilon_k^2.
$$

By solving this inequality, we obtain that

$$
\|\,\widehat{\mathcal{Q}}_{22}^{1/2}\,(\widetilde{v}^k-\overline{v}^k)\,\| \;\leq\; (\,1+\sqrt{2}\,)\,\varepsilon_k.
$$

(c) From parts (a) and (b), we have that

$$
\|\,\overline{w}^k-\widetilde{w}^k\,\|_{\mathcal{H}}^2 \;=\; \|\,\overline{u}^k-\widetilde{u}^k\,\|_{\mathcal{D}_1}^2+\|\,\overline{v}^k-\widetilde{v}^k\,\|_{\widehat{\mathcal{Q}}_{22}}^2 \;\leq\; 7\varepsilon_k^2,
$$

This completes the proof of this lemma. □

We are now ready to present the main theorem of this section on the $O(1/k^2)$ complexity of the imABCD algorithm.

**Theorem 3.** Suppose that Assumption 1 holds and the solution set $\Omega$ of the problem (3) is nonempty. Let $w^* \in \Omega$. Assume that $\sum_{i=1}^{\infty} i\,\varepsilon_i < \infty$. Then the sequence $\{\widetilde{w}^k\} \triangleq \{(\widetilde{u}^k, \widetilde{v}^k)\}$ generated by the imABCD algorithm satisfies that

$$\theta(\widetilde{w}^k) - \theta(w^*) \,\leq\, \frac{2\,\|\,\widetilde{w}^0 - w^*\,\|_{\mathcal{H}}^2 + c_0}{(k+1)^2}, \quad \forall\, k \geq 1,$$

where $c_0$ is a positive scalar (independent of $k$).

*Proof.* By applying Lemma 2 for $k = 1$ and $w = \overline{w}^1$, we get

$$\begin{aligned} 2\,[\,\theta(\overline{w}^1) - \theta(w^*)\,] \,&\leq\, \|\,w^1 - w^*\|_{\mathcal{H}}^2 - \|\,\overline{w}^1 - w^*\,\|_{\mathcal{H}}^2 \\ = \,\|\,\widetilde{w}^0 - w^*\|_{\mathcal{H}}^2 - \|\,\overline{w}^1 - w^*\|_{\mathcal{H}}^2 \,&=\, \|\,\overline{w}^0 - w^*\|_{\mathcal{H}}^2 - \|\,\overline{w}^1 - w^*\|_{\mathcal{H}}^2\,. \end{aligned}$$

For any $k \geq 2$, since $\overline{w}^k = (\overline{u}^k, \overline{v}^k)$ exactly solves the subproblem (7), we may take $\widetilde{w}^k = \overline{w}^k$ and $w = \dfrac{(t_k - 1)\overline{w}^{k-1} + w^*}{t_k}$ in Lemma 2 to obtain the following inequality:

$$\theta(w) - \theta(\overline{w}^k) \,\geq\, \frac{1}{2}\left\|\frac{(t_k-1)\overline{w}^{k-1} + w^*}{t_k} - \overline{w}^k\right\|_{\mathcal{H}}^2 - \frac{1}{2}\left\|\frac{(t_k-1)\overline{w}^{k-1} + w^*}{t_k} - w^k\right\|_{\mathcal{H}}^2.$$

By applying Lemma 1(c) with $w^+ = \overline{w}^{k-1}$ and $w' = w^*$, we get

$$t_k^2\,[\,\theta(w) - \theta(w^*)\,] \,\leq\, t_{k-1}^2\,[\,\theta(\overline{w}^{k-1}) - \theta(w^*)\,].$$

By combining the above two inequalities and noting that $t_k w^k = t_k\widetilde{w}^{k-1} + (t_{k-1} - 1)(\widetilde{w}^{k-1} - \widetilde{w}^{k-2})$, we have for $k \geq 2$,

$$\begin{aligned} &2\,t_k^2\left[\,\theta(\overline{w}^k) - \theta(w^*)\,\right] - 2\,t_{k-1}^2\left[\,\theta(\overline{w}^{k-1}) - \theta(w^*)\,\right] \\ \leq\;& \|\,t_{k-1}\widetilde{w}^{k-1} - w^* - (t_{k-1}-1)\widetilde{w}^{k-2} - (t_k - 1)(\overline{w}^{k-1} - \widetilde{w}^{k-1})\,\|_{\mathcal{H}}^2 - \|\,t_k\overline{w}^k - w^* - (t_k-1)\overline{w}^{k-1}\,\|_{\mathcal{H}}^2 \\ =\;& \|\,\lambda^{k-1}\,\|_{\mathcal{H}}^2 - 2\,\left\langle\,\mathcal{H}\lambda^{k-1}\,,\,(t_{k-1} + t_k - 1)(\overline{w}^{k-1} - \widetilde{w}^{k-1}) - (t_{k-1}-1)(\overline{w}^{k-2} - \widetilde{w}^{k-2})\,\right\rangle \\ &+ \|\,(t_{k-1} + t_k - 1)(\overline{w}^{k-1} - \widetilde{w}^{k-1}) - (t_{k-1}-1)(\overline{w}^{k-2} - \widetilde{w}^{k-2})\|_{\mathcal{H}}^2 \,-\, \|\,\lambda^k\,\|_{\mathcal{H}}^2\,, \end{aligned} \tag{16}$$

where $\lambda^k \,\triangleq\, t_k\overline{w}^k - w^* - (t_k - 1)\overline{w}^{k-1} \,=\, t_k(\overline{w}^k - w^*) - (t_k - 1)(\overline{w}^{k-1} - w^*)$. By Lemma 1(b), Lemma 4 and the nonincreasing property of $\{\varepsilon_k\}$, we derive that for all $k \geq 3$,

$$\begin{aligned} &\|\,(t_{k-1} + t_k - 1)(\overline{w}^{k-1} - \widetilde{w}^{k-1}) - (t_{k-1}-1)(\overline{w}^{k-2} - \widetilde{w}^{k-2})\,\|_{\mathcal{H}} \\ \leq\;& (t_{k-1} + t_k - 1)\,\|\,\overline{w}^{k-1} - \widetilde{w}^{k-1}\,\|_{\mathcal{H}} + (t_{k-1} - 1)\,\|\,\overline{w}^{k-2} - \widetilde{w}^{k-2}\,\|_{\mathcal{H}} \\ \leq\;& c_1(k-1)\,\varepsilon_{k-2}, \end{aligned}$$

where $c_1 = 5$. Note that in deriving the above inequality, we have used Lemma 4 (c). For $k = 2$, we also have that

$$\begin{aligned} &\|\,(t_{k-1} + t_k - 1)(\overline{w}^{k-1} - \widetilde{w}^{k-1}) - (t_{k-1}-1)(\overline{w}^{k-2} - \widetilde{w}^{k-2})\,\|_{\mathcal{H}} \,=\, t_2\|\,(\overline{w}^1 - \widetilde{w}^1)\|_{\mathcal{H}} \\ \leq\;& c_1(k-1)\varepsilon_{k-2}, \end{aligned}$$

where we set $\varepsilon_0 = \varepsilon_1$. It follows from (16) that for $k \geq 2$, we have

$$
\begin{aligned}
& 2\,t_k^2\,\left[\theta(\overline{w}^k) - \theta(w^*)\right] + \|\,\lambda^k\,\|_{\mathcal{H}}^2 \\
\leq\;& 2\,t_{k-1}^2\,\left[\theta(\overline{w}^{k-1}) - \theta(w^*)\right] + \|\,\lambda^{k-1}\,\|_{\mathcal{H}}^2 + 2c_1(k-1)\varepsilon_{k-2}\|\,\lambda^{k-1}\,\|_{\mathcal{H}} + c_1^2\,(k-1)^2\,\varepsilon_{k-2}^2 \\
\leq\;& \cdots \\
\leq\;& 2\,t_1^2\,\left[\theta(\overline{w}^1) - \theta(w^*)\right] + \|\,\lambda^1\,\|_{\mathcal{H}}^2 + 2c_1\sum_{i=1}^{k-1} i\,\varepsilon_{i-1}\,\|\,\lambda^i\,\|_{\mathcal{H}} + c_1^2\sum_{i=1}^{k-1} i^2\,\varepsilon_{i-1}^2 \\
\leq\;& \|\,\overline{w}^1 - w^*\,\|_{\mathcal{H}}^2 + 2c_1\sum_{i=1}^{k-1} i\,\varepsilon_{i-1}\,\|\,\lambda^i\,\|_{\mathcal{H}} + c_1^2\sum_{i=1}^{k-1} i^2\,\varepsilon_{i-1}^2\,.
\end{aligned}
$$

Notice that by Lemma 2, $\|\,\overline{w}^1 - w^*\,\|_{\mathcal{H}}^2 \leq \|\,w^1 - w^*\,\|_{\mathcal{H}}^2 = \|\,\overline{w}^0 - w^*\,\|_{\mathcal{H}}^2$. Next we show that the above inequality implies the boundedness of the sequence $\{\|\lambda^k\|_{\mathcal{H}}\}$. If $\|\lambda^k\|_{\mathcal{H}} \leq 1$ for all $k \geq 1$, then we are done. Otherwise, for any given sufficiently large positive integer $m$, we have that

$$
\|\,\lambda^{k_m}\,\|_{\mathcal{H}} \;:=\; \max\{\|\,\lambda^i\,\|_{\mathcal{H}} \mid 1 \leq i \leq m\} \geq 1.
$$

Thus, for any $1 \leq k \leq m$, we have that

$$
\begin{aligned}
\|\,\lambda^k\|_{\mathcal{H}} \;\leq\; \|\,\lambda^{k_m}\|_{\mathcal{H}} \;&\leq\; \frac{1}{\|\,\lambda^{k_m}\|_{\mathcal{H}}}\left(\|\,\overline{w}^0 - w^*\,\|_{\mathcal{H}}^2 + 2c_1\sum_{i=1}^{k_m-1} i\,\varepsilon_{i-1}\|\,\lambda^i\,\|_{\mathcal{H}} + c_1^2\sum_{i=1}^{k_m-1} i^2\,\varepsilon_{i-1}^2\right) \\
&\leq \|\,\overline{w}^0 - w^*\,\|_{\mathcal{H}}^2 + 2c_1\sum_{i=1}^{k_m-1} i\,\varepsilon_{i-1}\,\frac{\|\,\lambda^i\,\|_{\mathcal{H}}}{\|\,\lambda^{k_m}\,\|_{\mathcal{H}}} + c_1^2\sum_{i=1}^{k_m-1} i^2\,\varepsilon_{i-1}^2 \\
&\leq \|\,\overline{w}^0 - w^*\,\|_{\mathcal{H}}^2 + 2c_1\sum_{i=1}^{k_m-1} i\,\varepsilon_{i-1} + c_1^2\sum_{i=1}^{k_m-1} i^2\,\varepsilon_{i-1}^2 \\
&\leq \|\,\overline{w}^0 - w^*\,\|_{\mathcal{H}}^2 + 2c_1\sum_{i=1}^{\infty} i\,\varepsilon_{i-1} + c_1^2\sum_{i=1}^{\infty} i^2\,\varepsilon_{i-1}^2,
\end{aligned}
$$

where the third inequality follows from the definition of $\|\lambda^{k_m}\|_{\mathcal{H}}$. Thus by letting $m \to \infty$, we get

$$
\|\,\lambda^k\|_{\mathcal{H}} \;\leq\; c_2 \triangleq \max\left\{1,\; \|\,\overline{w}^0 - w^*\,\|_{\mathcal{H}}^2 + 2c_1\sum_{i=1}^{\infty} i\,\varepsilon_{i-1} + c_1^2\sum_{i=1}^{\infty} i^2\,\varepsilon_{i-1}^2\right\}, \quad \forall\, k \geq 1.
$$

To estimate the bound for $\|\,\overline{w}^{k+1} - w^*\|_{\mathcal{H}}$, we set $w = w^*$ and $\widetilde{w}^k = \overline{w}^k$ in Lemma 3 and deduce that

$$
\begin{aligned}
& t_{k+1}\,\|\,\overline{w}^{k+1} - w^*\,\|_{\mathcal{H}} \;\leq\; t_{k+1}\,\|\,w^{k+1} - w^*\,\|_{\mathcal{H}} \\
&= \left\|\,(t_{k+1} + t_k - 1)\widetilde{w}^k - (t_k - 1)\widetilde{w}^{k-1} - t_{k+1}w^*\,\right\|_{\mathcal{H}} \\
&\leq (t_{k+1} - 1)\left\|\,(\overline{w}^k - w^*)\,\right\|_{\mathcal{H}} + \left\|\,t_k\overline{w}^k - (t_k - 1)\overline{w}^{k-1} - w^*\,\right\|_{\mathcal{H}} \\
&\quad + (t_{k+1} + t_k - 1)\left\|\,(\widetilde{w}^k - \overline{w}^k)\,\right\|_{\mathcal{H}} + (t_k - 1)\left\|\,(\widetilde{w}^{k-1} - \overline{w}^{k-1})\,\right\|_{\mathcal{H}} \\
&\leq \frac{t_k^2}{t_{k+1}}\left\|\,(\overline{w}^k - w^*)\,\right\|_{\mathcal{H}} + c_2 + c_1\,\varepsilon_{k-1}\,t_{k+1}\,,
\end{aligned}
$$

where the last inequality is obtained by Lemma 1(a). It follows that

$$\begin{aligned}
\frac{t_k^2}{t_{k+1}^2}\left\| (\overline{w}^k - w^*)\right\|_{\mathcal{H}} &\le \frac{t_k^2}{t_{k+1}^2}\left(\frac{t_{k-1}^2}{t_k^2}\left\| (\overline{w}^{k-1} - w^*)\right\|_{\mathcal{H}} + \frac{c_2}{t_k} + c_1\varepsilon_{k-2}\right) \\
\frac{t_k^2}{t_{k+1}^2}\frac{t_{k-1}^2}{t_k^2}\left\| (\overline{w}^{k-1} - w^*)\right\|_{\mathcal{H}} &\le \frac{t_k^2}{t_{k+1}^2}\frac{t_{k-1}^2}{t_k^2}\left(\frac{t_{k-2}^2}{t_{k-1}^2}\left\| (\overline{w}^{k-2} - w^*)\right\|_{\mathcal{H}} + \frac{c_2}{t_{k-1}} + c_1\varepsilon_{k-3}\right) \\
\vdots \quad &\le \quad \vdots \\
\frac{t_k^2}{t_{k+1}^2}\frac{t_{k-1}^2}{t_k^2}\cdots\frac{t_2^2}{t_3^2}\left\| (\overline{w}^2 - w^*)\right\|_{\mathcal{H}} &\le \frac{t_k^2}{t_{k+1}^2}\frac{t_{k-1}^2}{t_k^2}\cdots\frac{t_2^2}{t_3^2}\left(\frac{t_1^2}{t_2^2}\left\| (\overline{w}^1 - w^*)\right\|_{\mathcal{H}} + \frac{c_2}{t_2} + c_1\varepsilon_0\right).
\end{aligned}$$

Summing up the above inequalities, we obtain

$$\left\| (\overline{w}^{k+1} - w^*)\right\|_{\mathcal{H}} \le \frac{t_1^2}{t_{k+1}^2}\left\| (\overline{w}^1 - w^*)\right\|_{\mathcal{H}} + c_2\sum_{i=1}^{k}\frac{t_{i+1}}{t_{k+1}^2} + c_1\sum_{i=1}^{k}\varepsilon_{i-1}. \tag{17}$$

By Lemma 1(b), we have

$$\sum_{i=1}^{k}\frac{t_{i+1}}{t_{k+1}^2} \le \frac{(3+k)\,k}{2\,\left(\frac{1}{2}k+1\right)^2} \le 2, \quad \forall\, k \ge 1.$$

Therefore, the inequality (17) implies that for all $k \ge 1$,

$$\begin{aligned}
\left\| (\overline{w}^{k+1} - w^*)\right\|_{\mathcal{H}} &\le \frac{4}{(k+2)^2}\left\| (\overline{w}^1 - w^*)\right\|_{\mathcal{H}} + 2c_2 + c_1\sum_{i=1}^{\infty}\varepsilon_{i-1} \\
&\le c_3 \triangleq \frac{4}{9}\left\| (\overline{w}^1 - w^*)\right\|_{\mathcal{H}} + 2c_2 + c_1\sum_{i=1}^{\infty}\varepsilon_{i-1}.
\end{aligned} \tag{18}$$

Note that we also have $\left\| (\overline{w}^1 - w^*)\right\|_{\mathcal{H}} \le c_3$.

The next step is to prove the boundedness of the term $\left\| t_k\widetilde{u}^k - u^* - (t_k - 1)\widetilde{u}^{k-1}\right\|_{\widehat{Q}_{11}}$. Before that, we need to first bound $\left\|\overline{u}^k - u^*\right\|_{\widehat{\mathcal{Q}}_{11}}^2$. By using Lemma 4, we have that for $k \ge 2$,

$$\begin{aligned}
& t_k\left\|\overline{u}^k - u^*\right\|_{\widehat{\mathcal{Q}}_{11}} \le t_k\left\|w^k - w^*\right\|_{\mathcal{H}} = \|(t_{k-1} + t_k - 1)(\widetilde{w}^{k-1} - w^*) - (t_{k-1} - 1)(\widetilde{w}^{k-2} - w^*)\|_{\mathcal{H}} \\
&\le (t_k - 1)\|(\widetilde{w}^{k-1} - w^*)\|_{\mathcal{H}} + \|t_{k-1}(\widetilde{w}^{k-1} - w^*) - (t_{k-1} - 1)(\widetilde{w}^{k-2} - w^*)\|_{\mathcal{H}} \\
&\le (t_k - 1)\left[\|\widetilde{w}^{k-1} - \overline{w}^{k-1}\|_{\mathcal{H}} + \|\overline{w}^{k-1} - w^*\|_{\mathcal{H}}\right] + \|t_{k-1}(\widetilde{w}^{k-1} - \overline{w}^{k-1}) - (t_{k-1} - 1)(\widetilde{w}^{k-2} - \overline{w}^{k-2})\|_{\mathcal{H}} \\
&\quad + \|t_{k-1}(\overline{w}^{k-1} - w^*) - (t_{k-1} - 1)(\overline{w}^{k-2} - w^*)\|_{\mathcal{H}} \\
&\le (t_k - 1)[\sqrt{7}\varepsilon_{k-1} + c_3\,] + (2t_k - 1)\sqrt{7}\varepsilon_{k-2} + \|\gamma^{k-1}\|_{\mathcal{H}} \\
&\le t_k[\,8\varepsilon_1 + c_2 + c_3].
\end{aligned}$$

For $k = 1$, we also have that $\|\overline{u}^k - u^*\|_{\widehat{\mathcal{Q}}_{11}} \le \|\overline{w}^k - w^*\|_{\widehat{\mathcal{H}}} \le c_3$.

Now we have for all $k \geq 2$,

$$
\begin{aligned}
& \| t_k \widetilde{u}^k - u^* - (t_k - 1)\, \widetilde{u}^{k-1} \|_{\widehat{\mathcal{Q}}_{11}} \\
\leq\ & t_k \Big( \| \overline{u}^k - \widetilde{u}^k \|_{\widehat{\mathcal{Q}}_{11}} + \| \overline{u}^k - u^* \|_{\widehat{\mathcal{Q}}_{11}} \Big) + (t_k - 1) \Big( \| \overline{u}^{k-1} - u^* \|_{\widehat{\mathcal{Q}}_{11}} + \| \overline{u}^{k-1} - \widetilde{u}^{k-1} \|_{\widehat{\mathcal{Q}}_{11}} \Big) \\
\leq\ & (2t_k - 1)\, [\, 8\, \varepsilon_1 + c_2 + c_3 \,].
\end{aligned}
$$

In addition, for $k = 1$, we have

$$
\begin{aligned}
& \| t_1 \widetilde{u}^1 - u^* - (t_1 - 1) \widetilde{u}^0 \|_{\widehat{Q}_{11}} \ = \ \| \widetilde{u}^1 - u^* \|_{\widehat{Q}_{11}} \ \leq \ \| \overline{u}^1 - u^* \|_{\widehat{Q}_{11}} + \| \widetilde{u}^1 - \overline{u}^1 \|_{\widehat{Q}_{11}} \\
& \quad \leq \ \| w^1 - w^* \|_{\mathcal{H}} + \varepsilon_1 \ \leq \ c_2 + \varepsilon_1 \ \leq \ (2t_k - 1)\, [\, 8\, \varepsilon_1 + c_2 + c_3 \,].
\end{aligned}
$$

On the other hand,

$$
\begin{aligned}
& \| t_k \widetilde{v}^k - v^* - (t_k - 1)\, \widetilde{v}^{k-1} \|_{\widehat{\mathcal{Q}}_{22}} \ = \ \| t_k (\widetilde{v}^k - v^*) - (t_k - 1)\, (\widetilde{v}^{k-1} - v^*) \|_{\widehat{\mathcal{Q}}_{22}} \\
\leq\ & \| t_k (\overline{v}^k - v^*) - (t_k - 1) (\overline{v}^{k-1} - v^*) \|_{\widehat{\mathcal{Q}}_{22}} + \| t_k (\widetilde{v}^k - \overline{v}^k) - (t_k - 1)\, (\widetilde{v}^{k-1} - \overline{v}^{k-1}) \|_{\widehat{\mathcal{Q}}_{22}} \\
\leq\ & \| \lambda^k \|_{\mathcal{H}} + t_k \| (\widetilde{v}^k - \overline{v}^k) \|_{\widehat{\mathcal{Q}}_{22}} + (t_k - 1) \| (\widetilde{v}^{k-1} - \overline{v}^{k-1}) \|_{\widehat{\mathcal{Q}}_{22}} \\
\leq\ & c_2 + (2t_k - 1)(1 + \sqrt{2}) \varepsilon_{k-1} \ \leq \ (2t_k - 1)[c_2 + 3\varepsilon_1].
\end{aligned}
$$

Finally, by applying Lemma 3 at $w = \dfrac{(t_k - 1)\widetilde{w}^{k-1} + w^*}{t_k}$ and using Lemma 1(c), we see that

$$
\begin{aligned}
& t_k^2 \big[ \theta(\widetilde{w}^k) - \theta(w^*) \big] + \frac{1}{2} \| t_k \widetilde{w}^k - w^* - (t_k - 1) \widetilde{w}^{k-1} \|_{\mathcal{H}}^2 \\
\leq\ & t_{k-1}^2 \big[ \theta(\widetilde{w}^{k-1}) - \theta(w^*) \big] + \frac{1}{2} \| t_{k-1} \widetilde{w}^{k-1} - w^* - (t_{k-1} - 1) \widetilde{w}^{k-2} \|_{\mathcal{H}}^2 \\
& + \varepsilon_k \| t_k \widetilde{w}^k - w^* - (t_k - 1) \widetilde{w}^{k-1} \|_{\mathrm{Diag}(\widehat{Q}_{11},\, \widehat{Q}_{22})} \\
\leq\ & \cdots \\
\leq\ & t_1^2 \big[ \theta(\widetilde{w}^1) - \theta(w^*) \big] + \frac{1}{2} \| t_1 \widetilde{w}^1 - w^* - (t_1 - 1) \widetilde{w}^0 \|_{\mathcal{H}}^2 + \sum_{i=2}^{k} \varepsilon_i \| t_i \widetilde{w}^i - w^* - (t_i - 1) \widetilde{w}^{i-1} \|_{\mathrm{Diag}(\widehat{Q}_{11},\, \widehat{Q}_{22})} \\
\leq\ & \frac{1}{2} \| \widetilde{w}^0 - w^* \|_{\mathcal{H}}^2 + \sum_{i=1}^{k} \varepsilon_i \| t_i \widetilde{w}^i - w^* - (t_i - 1) \widetilde{w}^{i-1} \|_{\mathrm{Diag}(\widehat{Q}_{11},\, \widehat{Q}_{22})} \\
\leq\ & \frac{1}{2} \| \widetilde{w}^0 - w^* \|_{\mathcal{H}}^2 + \sum_{i=1}^{k} (2t_i - 1)\, [\, 11\varepsilon_1 + 2c_2 + c_3 \,]\, \varepsilon_i \\
\leq\ & \frac{1}{2} \| \widetilde{w}^0 - w^* \|_{\mathcal{H}}^2 + \frac{1}{4} c_0,
\end{aligned}
$$

where

$$
c_0 \triangleq 4 \sum_{i=1}^{\infty} (2t_i - 1)\, [\, 11\varepsilon_1 + 2c_2 + c_3 \,]\, \varepsilon_i
$$

is a finite value since $\sum_{i=1}^{\infty} t_i\, \varepsilon_i < \infty$. Since $t_k \geq \dfrac{k+1}{2}$, we complete the proof of this theorem. □

# 4 Solving the Best Approximation Problem (2)

In this section, we discuss an application of the imABCD framework to solve the dual of the best approximation problem (2). The study of the best approximation problems with only equalities constraints dates back three decades [19, 2, 9, 10, 11]. The best approximation problem with the positive semidefinite cone constraint was studied recently, see, e.g., [25, 31, 6, 16, 39].

The dual of (2) is given by

$$\underset{y,S,z,Z}{\text{minimize}} \quad \frac{1}{2}\|\mathcal{A}^*y + S + \mathcal{B}^*z + Z + G\,\|^2 - \langle b, y\rangle - \langle d, z\rangle + \delta_{\mathbb{S}^n_+}(S) + \delta_{\geq 0}(z) + \delta_{\geq 0}(Z), \tag{19}$$

where $\delta_C(\cdot)$ denotes the indicator function of a given set $C$, i.e., $\delta_C(x) = 0$ if $x \in C$, and $\delta_C(x) = \infty$ otherwise. The notation $\delta_{\geq 0}(\cdot)$ is used to denote the indicator function over a nonnegative orthant. To implement the two-block imABCD algorithm, we take $(y, S)$ as one block, and $(z, Z)$ as the other block. There are two important reasons for us to merge the four blocks in such a way. First, based on our previous experiences from developing several solvers with similar types of constraints, in particular for linear and least squares semidefinite programming in [48, 39], we find that, compared to the linear inequalities and the nonnegative cone constraints, the linear equalities and the semidefinite cone constraints are more challenging to be satisfied numerically. Putting the corresponding multipliers $y$ (for the linear equalities constraints) and $S$ (for the positive semidefinite cone constraint) in one group and solving them simultaneously by the inexact semismooth Newton method may help these two types of constraints to achieve a high accuracy simultaneously. Second, putting $S$ and $z$ (corresponding to the linear inequalities constraints) or $Z$ (corresponding to the entrywise nonnegative constraints) in one group often leads to a degenerate Newton system when the semismooth Newton method is applied to solve the resulting subproblem.

## 4.1 Solving the subproblems by the Newton-type methods

With fixed $(z, Z)$ and a properly defined matrix $G_1 \in \mathbb{S}^n$, the first workhorse of the imABCD for solving (2) is the following convex program:

$$\underset{y,S}{\text{minimize}} \quad \frac{1}{2}\,\|\mathcal{A}^*y + S + G_1\|^2 - \langle b, y\rangle \quad \text{subject to} \quad S \in \mathbb{S}^n_+,\ y \in \mathbb{R}^m.$$

Let $\Pi_{\mathbb{S}_n^+}(\cdot)$ denote the projection onto the cone $\mathbb{S}^n_+$. Since the optimal solution $(\bar{y}, \overline{S})$ always satisfies

$$\overline{S} \;=\; \Pi_{\mathbb{S}^n_+}(-G_1 - \mathcal{A}^*\bar{y}), \tag{20}$$

we may solve $y$ first via the following unconstrained minimization problem:

$$\underset{y}{\text{minimize}}\ \xi(y) \;\triangleq\; \frac{1}{2}\,\|\mathcal{A}^*y + G_1 + \Pi_{\mathbb{S}^n_+}(-G_1 - \mathcal{A}^*y)\|^2 - \langle b, y\rangle \;=\; \frac{1}{2}\,\|\Pi_{\mathbb{S}^n_+}(G_1 + \mathcal{A}^*y)\|^2 - \langle b, y\rangle \tag{21}$$

and then substitute the solution into (20) to obtain $\overline{S}$. Notice that $\xi$ is a continuously differentiable function and its gradient

$$\nabla\xi(y) \;=\; \mathcal{A}\Pi_{\mathbb{S}^n_+}(G_1 + \mathcal{A}^*y) - b$$

is strongly semismooth [38]. Therefore, the semismooth Newton-CG algorithm with line search, which is proven to converge globally and local superlinearly/quadratically [22, 32, 31, 49], can be applied to solve the above unconstrained problem (21). The details of this algorithm are given below.

---

**SNCG:** a **S**emismooth **N**ewton-**CG** method for solving (21)

---

**Initialization.** Given $\mu \in (0, 1/2)$, $\eta \in (0, 1)$, $\tau \in (0, 1]$ and $\rho \in (0, 1)$. Iterate the following steps for $j \geq 0$.

**Step 1.** Choose $V^j \in \partial \Pi_{\mathbb{S}^n_+}(G_1 + \mathcal{A}^* y^j)$. Solve the following linear system to find $d^j$ by the conjugate gradient (CG) method:

$$\mathcal{A} V^j \mathcal{A}^* \, d + \nabla \xi(y^j) = 0$$

such that $d^j$ satisfies the accuracy condition that $\| \, \mathcal{A} V^j \mathcal{A}^* \, d + \nabla \xi(y^j) \, \| \; \leq \; \min\{ \, \eta, \| \, \nabla \xi(y^j) \, \|^{1+\tau} \, \}$.

**Step 2.** (Line search) Set $\alpha_j \; = \; \rho^{m_j}$, where $m_j$ is the first nonnegative integer $m$ for which

$$\xi(y^j + \rho^m d_j) \; \leq \; \xi(y^j) + \mu \, \rho^m \, \langle \nabla \xi(y^j), d_j \rangle.$$

**Step 3.** Set $y^{j+1} \; = \; y^j + \alpha_j d^j$.

---

To solve the second subproblem involving $(z, Z)$ in (7), we need to deal with the program

$$\underset{z,Z}{\text{minimize}} \quad \frac{1}{2} \, \| \, \mathcal{B}^* z + Z + G_2 \, \|^2 - \langle d, z \rangle + \frac{c}{2} \, \| \, z - z_0 \, \|^2 \quad \text{subject to } z \geq 0, \quad Z \geq 0 \tag{22}$$

for some properly defined matrix $G_2 \in \mathbb{S}^n$. Notice that an additional proximal term $\frac{c}{2} \| z - z_0 \|^2$ for some scalar $c > 0$ is added to the objective function (which serves as a majorant of the original objective function) to make this problem strongly convex with respect to $z$. Similar to the first subproblem, the optimal solution $(\bar{z}, \overline{Z})$ satisfies that

$$\overline{Z} \; = \; \Pi_{\geq 0}(-\mathcal{B}^* \bar{z} - G_2).$$

We therefore need to solve the problem

$$\underset{z}{\text{minimize}} \quad \frac{1}{2} \, \| \, \Pi_{\geq 0}(\mathcal{B}^* z + G_2) \, \|^2 - \langle d, z \rangle + \frac{c}{2} \, \| \, z - z_0 \, \|^2 \quad \text{subject to } z \geq 0. \tag{23}$$

Different from (21), the above problem is a constrained $\text{SC}^1$ problem (i.e., the objective function is continuously differentiable with a semismooth gradient). Though we may still apply a globally convergent semismooth Newton algorithm as proposed in [30], however, a strictly convex quadratic programming problem has to be solved in each step, which itself may be challenging. In fact, the problem (23) is a special case of the general unconstrained nonsmooth convex program

$$\underset{x}{\text{minimize}} \; \zeta(x) + \psi(x), \tag{24}$$

where $\zeta : \mathbb{X} \to (-\infty, \infty)$ is a strongly convex and smooth function, and $\psi : \mathbb{X} \to (-\infty, +\infty]$ is a convex but possibly nonsmooth function. One can apply Nesterov's accelerated proximal gradient method (APG) [27], which converges globally and linearly, to solve such a strongly convex problem [37]. Alternatively, the solution of (24) can be obtained via the nonsmooth equation

$$F(x) \triangleq x - \mathrm{Prox}_{\psi}\,(\, x - \nabla\zeta(x)\,) = 0,$$

where $\mathrm{Prox}_{\psi}(x) \triangleq \mathrm{argmin}\left\{\psi(x') + \frac{1}{2}\,\|\, x' - x\,\|^2 \;\middle|\; x' \in \mathbb{X}\right\}$ denotes the proximal mapping of $\psi$ at $x$, see, e.g., [35, Definition 1.22]. Since the composition of semismooth functions is semismooth [15], it follows that $F$ is semismooth at $x$ if $\nabla\zeta$ is semismooth at $x$ and $\mathrm{Prox}_{\psi}(\cdot)$ is semismooth at $x - \nabla\zeta(x)$. We may then apply the semismooth Newton-CG method locally to solve the above nonsmooth equation for a faster convergence rate. Therefore, a convergent and efficient way to solve (23) may be a hybrid of the APG algorithm and the semismooth Newton-CG method. We present this algorithm below, where the positive scalar $L_{\zeta}$ denotes the Lipschitz constant of $\nabla\zeta$.

---

**APG-SNCG**: A hybrid of the **APG** algorithm and the **SNCG** method for solving (24)

---

Choose an initial point $x^1 \in \mathbb{X}$, positive constants $\eta, \gamma \in (0,1)$, $\rho \in (0, 1/2)$, and a positive integer $m_0 > 0$. Iterate the following steps for $j \geq 0$.

**Step 1:** Select $V^j \in \partial F(x^j)$, the generalized Jacobian of $F$ at $x^j$, and apply the CG method to find an approximate solution $d^j$ to

$$V^j d + F(x^j) = 0 \tag{25}$$

such that

$$R^j \triangleq V^j d^j + F(x^j) \quad \text{and} \quad \|\, R^j \,\| \;\leq\; \eta_j \,\|\, F(x^j)\,\|, \tag{26}$$

where $\eta_j \triangleq \min\{\,\eta, \|\, F(x^j)\,\|\,\}$. If (26) is achievable, go to Step 2. Otherwise, go to Step 3.

**Step 2:** Let $m_j \leq m_0$ be the smallest nonnegative integer $m$ such that

$$\|\, F(x^j + \rho^m d^j)\,\| \;\leq\; \gamma \,\|\, F(x^j)\,\|.$$

If the above inequality is achievable, set $t_j = \rho^{m_j}$ and $x^{j+1} = x^j + t_j\, d^j$. Replace $j$ by $j+1$ and go to Step 1. Otherwise (i.e., the above inequality fails for all $m \leq m_0$) go to step 3.

**Step 3:** Set $x^{j_1} = \tilde{x}^{j_0} = x^j$, $\beta_{j_1} = 1$ and $i = 1$, compute

$$\begin{cases} \tilde{x}^{j_i} \;=\; \mathrm{Prox}_{\psi/L_{\zeta}}(x^{j_i} - \nabla\zeta(x^{j_i})/L_{\zeta}), \\ \beta_{j_{i+1}} \;=\; \dfrac{1}{2}\left(1 + \sqrt{1 + 4\beta_{j_i}^2}\right), \\ x^{j_{i+1}} \;=\; \tilde{x}^{j_i} + \dfrac{\beta_{j_i} - 1}{\beta_{j_{i+1}}}(\tilde{x}^{j_i} - \tilde{x}^{j_{i-1}}). \end{cases}$$

If $\|\, F(x^{j_{i+1}})\,\| \;\leq\; \gamma \,\|\, F(x^j)\,\|$, set $x^{j+1} = x^{j_{i+1}}$. Replace $j$ by $j+1$ and go to Step 1. Otherwise, set $i = i + 1$ and continue the above iteration.

---

Since $\zeta$ is strongly convex, the sequence $\{x^{j_i}\}$ generated by the APG algorithm always converges to the unique optimal solution $x^*$ of the problem (24), which further implies that $F(x^{j_i}) \to 0$ by the continuity of the proximal mapping. Therefore, the APG algorithm can be viewed as a safeguard for the global convergence of $\{x^j\}$ in the above framework. We make two further remarks of this algorithm below.

**Remark 4.** It is known from Rademacher's Theorem that the Lipschitz continuous function $F$ is differentiable almost everywhere. Assume that (26) is achievable at a differentiable point $x^j$ and $\|F(x^j)\| \neq 0$, then $\|F(x)\|^2$ is differentiable at $x^j$ and

$$\begin{aligned} \| F(x^j + td^j) \|^2 &= \| F(x^j) + t\,[\, R^j - F(x^j)\,] + o(t) \|^2 \\ &= \| F(x^j) \|^2 + t \,\langle\, F(x^j)\,,\, R^j - F(x^j) \,\rangle + o(t) \\ &\leq \| F(x^j) \|^2 + t\,(\eta_j - 1) \,\| F(x^j) \|^2 + o(t). \end{aligned}$$

Since $\eta_j \leq \eta < 1$, we have $\| F(x^j + t\,d^j) \| < \| F(x^j) \|$ for $t$ sufficiently small such that $d^j$ is a descent direction of $\| F(x) \|$ at $x^j$. This yields that the direction obtained by (25) is a descent direction of $\| F(x) \|$ at $x^j$ with probability 1.

**Remark 5.** The equation (25) may not be a symmetric linear system. If this occurs, one may use the BiCGStab iterative solver (e.g., van der Vorst [44]) to solve the corresponding equation.

## 4.2 Decomposing (22) into smaller decoupled problems

For some best approximation problems, the number of inequalities in $\mathcal{B}X \geq d$ may be ultra large, and that can make the subproblem (22) extremely expensive to solve. Fortunately, by the design of imABCD, one can add an appropriate proximal term in (7) to make the subproblem involving $(z, Z)$ easier to solve. In particular, by decomposing (22) into smaller decoupled problems.

A practical way to achieve the decomposition of (22) is to add a proximal term of the form $\frac{1}{2}\| \begin{pmatrix} z \\ Z \end{pmatrix} - \begin{pmatrix} z_0 \\ Z_0 \end{pmatrix} \|^2_{\widehat{\mathcal{Q}}}$, where the positive semidefinite linear operator $\widehat{\mathcal{Q}}$ is constructed based on dividing the operator $\mathcal{B}$ and the dual variable $z$ into $q \geq 1$ parts as

$$\begin{cases} \mathcal{B}X \equiv \begin{pmatrix} \mathcal{B}_1 X \\ \vdots \\ \mathcal{B}_q X \end{pmatrix} \text{ with } \mathcal{B}_i : \mathbb{S}^n \to \mathbb{R}^{m_i}, \quad X \in \mathbb{S}^n, \\ \mathcal{B}^* z \equiv (\mathcal{B}_1^* z_1, \mathcal{B}_2^* z_2, \ldots, \mathcal{B}_q^* z_q), \quad z \equiv (z_1, z_2, \ldots z_q) \in \mathbb{R}^{m_{i_1}} \times \mathbb{R}^{m_{i_2}} \times \ldots \times \mathbb{R}^{m_{i_q}} = \mathbb{R}^{m_I}. \end{cases}$$

By observing that for any given matrix $X \in \mathbb{R}^{m\times n}$,

$$\begin{pmatrix} & X \\ X^T & \end{pmatrix} \preceq \begin{pmatrix} (XX^T)^{\frac{1}{2}} & \\ & (X^TX)^{\frac{1}{2}} \end{pmatrix},$$

we may derive that

$$\mathcal{B}\mathcal{B}^* \preceq \mathcal{M} \triangleq \operatorname{Diag}(\mathcal{M}_1, \ldots, \mathcal{M}_q),$$

where

$$\mathcal{M}_i \triangleq \mathcal{B}_i\mathcal{B}_i^* + \sum_{j=1,\ldots,q,\, j\neq i} (\mathcal{B}_i\mathcal{B}_j^*\mathcal{B}_j\mathcal{B}_i^*)^{1/2}, \quad i = 1, \ldots q.$$

It therefore follows that

$$\begin{pmatrix} 2\mathcal{M} & \\ & 2I \end{pmatrix} \succeq \begin{pmatrix} 2\mathcal{B}\mathcal{B}^* & \\ & 2I \end{pmatrix} \succeq \begin{pmatrix} \mathcal{B} \\ I \end{pmatrix} \begin{pmatrix} \mathcal{B}^* & I \end{pmatrix}.$$

By choosing $\widehat{\mathcal{Q}} = \operatorname{diag}(2\mathcal{M}, 2I) - (\mathcal{B}^*,\ I)^*(\mathcal{B}^*,\ I) \succeq 0$, we can show that

$$
\begin{aligned}
&\frac{1}{2}\left\|(\mathcal{B}^*,\ I)\begin{pmatrix} z \\ Z\end{pmatrix} + G_2\right\|^2 - \langle d, z\rangle + \frac{c}{2}\|\, z - z_0 \,\|^2 + \frac{1}{2}\left\| \begin{pmatrix} z \\ Z\end{pmatrix} - \begin{pmatrix} z_0 \\ Z_0\end{pmatrix} \right\|^2_{\widehat{\mathcal{Q}}} \\
=\ \ & \big\langle z,\, \mathcal{M}z\big\rangle + \frac{c}{2}\|z\|^2 - \big\langle z,\, h\big\rangle + \big\langle Z,\, Z\big\rangle + \big\langle Z,\, G_2 - Z_0 + \mathcal{B}^* z_0\big\rangle + \kappa \\
=\ \ & \kappa + \big\langle Z,\, Z\big\rangle + \big\langle Z,\, G_2 - Z_0 + \mathcal{B}^* z_0\big\rangle + \sum_{i=1}^{q} \frac{1}{2}\big\langle (2\mathcal{M}_i + c\, I) z_i,\, z_i\big\rangle - \big\langle z_i,\, h_i\big\rangle,
\end{aligned}
$$

where $h = d + cz_0 + 2\mathcal{M}z_0 - \mathcal{B}(G_2 + Z_0 + \mathcal{B}^* z_0)$ and $\kappa = \frac{c}{2}\|z_0\|^2 + \frac{1}{2}\left\| \begin{pmatrix} z_0 \\ Z_0\end{pmatrix}\right\|^2_{\widehat{\mathcal{Q}}} + \frac{1}{2}\|G_2\|^2$. Thus, we can obtain the following decomposition for solving (22) with an additional proximal term:

$$
\begin{aligned}
&\underset{z,Z}{\text{minimize}} \left\{ \frac{1}{2}\|\, \mathcal{B}^* z + Z + G_2 \,\|^2 - \langle d, z\rangle + \frac{c}{2}\|\, z - z_0 \,\|^2 + \frac{1}{2}\left\| \begin{pmatrix} z \\ Z\end{pmatrix} - \begin{pmatrix} z_0 \\ Z_0\end{pmatrix} \right\|^2_{\widehat{\mathcal{Q}}} \ \Big|\ z \geq 0,\ Z \geq 0 \right\} \\
\Longleftrightarrow\ \ & \underset{Z}{\text{minimize}} \left\{ \big\langle Z,\, Z\big\rangle + \big\langle Z,\, G_2 - Z_0 + \mathcal{B}^* z_0\big\rangle \mid Z \geq 0 \right\} \\
& + \sum_{i=1}^{q} \underset{z_i}{\text{minimize}} \left\{ \frac{1}{2}\big\langle (2\mathcal{M}_i + cI) z_i,\, z_i\big\rangle - \big\langle z_i,\, h_i\big\rangle \ \mid\ z_i \geq 0 \right\}.
\end{aligned}
$$

Observe that the subproblem with respect to $Z$ can be solved analytically. For each $i$, the decoupled subproblem with respect to $z_i$ is simple convex quadratic program of the form (24) that can be solved by the APG-SNCG method.

# 5 Numerical Experiments

In this section, we test our imABCD algorithm on solving the dual problem (19). The equalities and inequalities constraints are generated from the doubly nonnegative relaxation of a binary integer nonconvex quadratic (ex-BIQ) programming that was considered in [39]:

$$
\begin{aligned}
\underset{x,Y,X}{\text{minimize}} \quad & \frac{1}{2}\langle Q, Y\rangle + \langle c, x\rangle \\
\text{subject to} \quad & \operatorname{Diag}(Y) = x, \quad \alpha = 1, \quad X = \begin{pmatrix} Y & x \\ x^T & \alpha \end{pmatrix}, \quad X \in \mathbb{S}^n_+, \quad X \geq 0, \\
& -Y_{ij} + x_i \geq 0, \quad -Y_{ij} + x_j \geq 0, \quad Y_{ij} - x_i - x_j \geq -1, \quad \forall\, i < j,\ j = 2, \ldots, n-1.
\end{aligned}
$$

The test data for $Q$ and $c$ are taken from Biq Mac Library maintained by Wiegele, which is available at http://biqmac.uni-klu.ac.at/biqmaclib.html.

We set $G = -\frac{1}{2}\begin{pmatrix} Q & c \\ c & 0 \end{pmatrix}$ in (2). Under Slater's condition, the KKT optimality conditions of the problem (2) are:

$$
\begin{cases}
X = G + \mathcal{A}^* y + \mathcal{B}^* z + S + Z, \quad \mathcal{A}X = b, \\
\mathcal{B}X - d = \Pi_{\geq 0}(\mathcal{B}X - d - z), \quad X = \Pi_{\mathbb{S}^n_+}(X - S), \quad X = \Pi_{\geq 0}(X - Z).
\end{cases}
$$

We measure the accuracy of an approximate dual solution $(y, z, S, Z)$ by the relative residual of the KKT system $\eta \triangleq \max\{\eta_1, \eta_2, \eta_3, \eta_4\}$, where, by letting $X = G + \mathcal{A}^*y + \mathcal{B}^*z + S + Z$,

$$\begin{cases} \eta_1 \triangleq \dfrac{\|\mathcal{A}X - b\|}{1 + \|b\|}, \quad \eta_2 \triangleq \dfrac{\|\mathcal{B}X - d - \Pi_{\geq 0}(\mathcal{B}X - d - z)\|}{1 + \|d\|}, \\ \eta_3 \triangleq \dfrac{\|X - \Pi_{\mathbb{S}^n_+}(X - S)\|}{1 + \|X\| + \|S\|}, \quad \eta_4 \triangleq \dfrac{\|X - \Pi_{\geq 0}(X - Z)\|}{1 + \|X\| + \|Z\|}. \end{cases}$$

We also report the duality gap defined by

$$\eta_{\rm g} \triangleq \frac{\mathrm{obj_p} - \mathrm{obj_d}}{1 + |\mathrm{obj_p}| + |\mathrm{obj_d}|},$$

where $\mathrm{obj_p} \triangleq \frac{1}{2}\|X - G\|^2$ and $\mathrm{obj_d} \triangleq -\frac{1}{2}\|\mathcal{A}^*y + \mathcal{B}^*z + S + Z + G\|^2 + \langle b, y\rangle + \langle d, z\rangle + \frac{1}{2}\|G\|^2$. We stop the algorithms if $\eta < \varepsilon$ for some prescribed accuracy $\varepsilon$.

In order to demonstrate the importance for incorporating the second order information when solving the subproblems, we compare our imABCD method with the two-block accelerated block coordinate gradient descent algorithm proposed by Chambolle and Pock in [7]. For a fair comparison, the two blocks are also taken as $(y, S)$ and $(z, Z)$. Let $\lambda_{\max}(\mathcal{BB}^*)$ be the largest eigenvalue of $\mathcal{BB}^*$. The algorithm adapted from [7] is given as follows:

---

**ABCGD:** An **A**ccelerated **B**lock **C**oordinate **G**radient **D**escent algorithm for solving (19).

---

**Initialization**. Choose an initial point $W^1 = \widetilde{W}^0 \in \mathbb{W}$. Let $t_1 = 1$.

**Step 1.** Compute $R^{k+\frac{1}{2}} = \mathcal{A}^*y^k + \mathcal{B}^*z^k + S^k + Z^k + G$ and
$$\begin{cases} \widetilde{y}^k = y^k - \frac{1}{2}(\mathcal{AA}^*)^{-1}(\mathcal{A}R^{k+\frac{1}{2}} - b), \\ \widetilde{S}^k = \Pi_{\mathbb{S}^n_+}\left(S^k - \frac{1}{2}R^{k+\frac{1}{2}}\right). \end{cases}$$

**Step 2.** Compute $R^k = \mathcal{A}^*\widetilde{y}^k + \mathcal{B}^*z^k + \widetilde{S}^k + Z^k + G$ and
$$\begin{cases} \widetilde{z}^k = \Pi_{\geq 0}\left(z^k - \dfrac{1}{2\lambda_{\max}(\mathcal{BB}^*)}(\mathcal{B}R^k - d)\right), \\ \widetilde{Z}^k = \Pi_{\geq 0}\left(Z^k - \frac{1}{2}R^k\right). \end{cases}$$

**Step 3.** Compute
$$\begin{cases} t_{k+1} = \frac{1}{2}\left(1 + \sqrt{1 + 4t_k^2}\right), \\ W^{k+1} = \widetilde{W}^k + \dfrac{t_k - 1}{t_{k+1}}(\widetilde{W}^k - \widetilde{W}^{k-1}). \end{cases}$$

---

Figure 1 shows the performance profile of the imABCD and ABCGD algorithms for some large scale ex-BIQ problems with $\varepsilon = 10^{-6}$. Recall that a point $(x, y)$ is on the performance profile curve of a method if and only if it can solve exactly $(100y)\%$ of all the tested problems at most $x$ times slower than any other methods. The detailed numerical results are presented in Table 1. The first four columns list the problem names, the dimension of the variable $y$ ($m_E$), $z$ ($m_I$) and the size of the matrix $G$ ($n_s$), respectively. The last several columns report the number of iterations, the relative residual $\eta$, the relative duality gap $\eta_{\rm gap}$, and the computation time in the format of "hours:minutes:seconds". One can see from the performance profile that the ABCGD algorithm requires at least 5 times the number of iterations taken by imABCD and is about 4 times slower than the imABCD algorithm. In particular, the ABCGD method cannot solve all the large scale bdq500 problems within 50000 iterations, whereas our imABCD can obtain satisfactory solutions after 6000 iterations. This indicates that even though the computational cost for each cycle of the imABCD method is larger than that of the ABCGD method, this cost is compensated by taking much fewer iterations. In fact, the Newton system is well-conditioned in this case such that it only takes only one or two CG iterations to compute a satisfactory Newton direction.

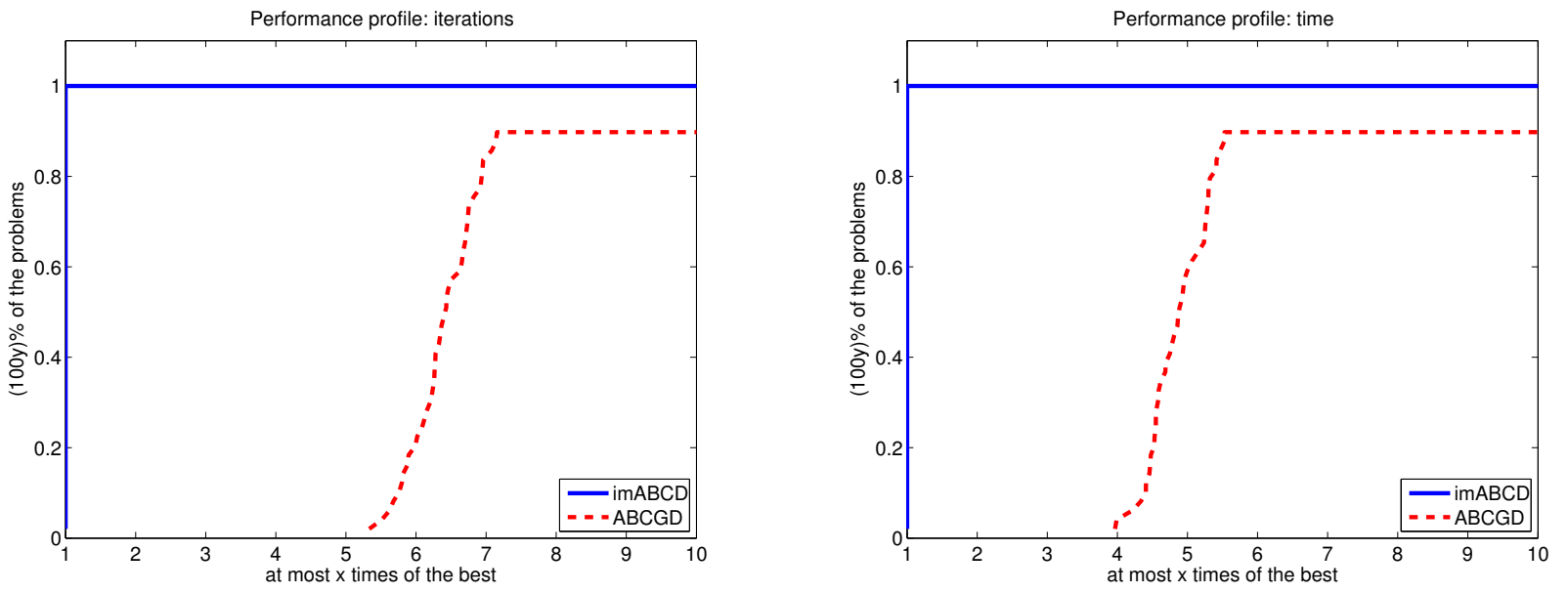

Figure 1: Performance profile of imABCD and ABCGD with $\varepsilon = 10^{-6}$

| | | | iterations | | $\eta$ | | $\eta_{\text{gap}}$ | | time | |
|---|---|---|---|---|---|---|---|---|---|---|
| problem | $m_E; m_I$ | $n_s$ | imabcd | abcgd | imabcd | abcgd | imabcd | abcgd | imabcd | abcgd |
| be100.1 | 101 ; 14850 | 101 | 5276 | 31048 | 9.9-7 | 9.9-7 | -2.3-7 | -7.7-8 | 45 | 3:08 |
| be100.2 | 101 ; 14850 | 101 | 5747 | 30645 | 9.9-7 | 9.9-7 | -2.0-7 | -9.4-8 | 47 | 3:08 |
| be100.3 | 101 ; 14850 | 101 | 5950 | 41172 | 9.9-7 | 9.9-7 | -3.7-7 | -8.7-8 | 53 | 4:02 |
| be100.4 | 101 ; 14850 | 101 | 5704 | 36684 | 9.9-7 | 9.9-7 | -2.7-7 | -7.2-8 | 49 | 3:48 |
| be100.5 | 101 ; 14850 | 101 | 5762 | 39956 | 9.9-7 | 9.7-7 | -3.1-7 | -7.6-8 | 51 | 4:05 |
| be100.6 | 101 ; 14850 | 101 | 5769 | 36134 | 9.9-7 | 9.9-7 | -1.8-7 | -7.9-8 | 48 | 3:40 |
| be100.7 | 101 ; 14850 | 101 | 4994 | 28087 | 9.9-7 | 9.9-7 | -3.3-7 | -8.4-8 | 40 | 2:53 |
| be100.8 | 101 ; 14850 | 101 | 5613 | 33772 | 9.9-7 | 9.9-7 | -1.5-7 | -5.3-8 | 45 | 3:32 |
| be100.9 | 101 ; 14850 | 101 | 5763 | 40048 | 9.9-7 | 9.9-7 | -2.7-7 | -8.1-8 | 51 | 4:10 |
| be100.10 | 101 ; 14850 | 101 | 5260 | 32010 | 9.9-7 | 9.9-7 | -3.7-7 | -8.8-8 | 43 | 3:46 |
| be120.3.1 | 121 ; 21420 | 121 | 4120 | 27781 | 9.9-7 | 9.9-7 | -2.2-7 | -5.6-8 | 41 | 3:42 |
| be120.3.2 | 121 ; 21420 | 121 | 4106 | 23809 | 9.9-7 | 9.9-7 | -2.3-7 | -6.3-8 | 40 | 3:00 |
| be120.3.3 | 121 ; 21420 | 121 | 3548 | 21867 | 9.9-7 | 9.9-7 | -9.8-8 | -6.6-8 | 35 | 2:53 |
| be120.3.4 | 121 ; 21420 | 121 | 4745 | 31783 | 9.9-7 | 9.9-7 | -3.1-7 | -6.7-8 | 47 | 4:08 |
| be120.3.5 | 121 ; 21420 | 121 | 5637 | 31076 | 9.9-7 | 9.9-7 | -4.9-8 | -7.1-8 | 58 | 3:51 |
| be120.3.6 | 121 ; 21420 | 121 | 3946 | 26558 | 9.9-7 | 9.9-7 | -1.8-7 | -4.8-8 | 39 | 3:26 |
| be120.3.7 | 121 ; 21420 | 121 | 4169 | 26176 | 9.9-7 | 9.9-7 | -2.7-7 | -6.5-8 | 41 | 3:37 |
| be120.3.8 | 121 ; 21420 | 121 | 3793 | 23796 | 9.9-7 | 9.9-7 | -1.6-7 | -3.8-8 | 35 | 3:11 |
| be120.3.9 | 121 ; 21420 | 121 | 4951 | 28518 | 9.9-7 | 9.9-7 | -2.0-7 | -5.2-8 | 52 | 3:58 |
| be120.3.10 | 121 ; 21420 | 121 | 4264 | 24803 | 9.9-7 | 9.9-7 | -3.4-7 | -5.3-8 | 42 | 3:06 |
| be120.8.1 | 121 ; 21420 | 121 | 5671 | 32200 | 9.9-7 | 9.9-7 | -3.5-7 | -8.1-8 | 58 | 4:26 |
| be120.8.2 | 121 ; 21420 | 121 | 5897 | 35336 | 9.9-7 | 9.9-7 | -3.1-7 | -7.1-8 | 1:02 | 4:39 |
| be120.8.3 | 121 ; 21420 | 121 | 5199 | 33259 | 9.9-7 | 9.9-7 | -4.9-7 | -9.2-8 | 52 | 4:35 |
| be120.8.4 | 121 ; 21420 | 121 | 6688 | 40964 | 9.9-7 | 9.9-7 | -4.0-7 | -8.0-8 | 1:12 | 5:26 |
| be120.8.5 | 121 ; 21420 | 121 | 5828 | 41263 | 9.9-7 | 9.9-7 | -3.6-7 | -6.3-8 | 1:03 | 5:40 |
| be120.8.6 | 121 ; 21420 | 121 | 4735 | 29524 | 9.9-7 | 9.9-7 | -4.8-7 | -7.8-8 | 47 | 3:52 |
| be120.8.7 | 121 ; 21420 | 121 | 4456 | 29722 | 9.9-7 | 9.9-7 | -3.8-7 | -6.5-8 | 44 | 3:36 |

| | | iterations | $\eta$ | $\eta_{\text{gap}}$ | time |
|---|---|---|---|---|---|
| problem | $m_E; m_I \mid n_s$ | imabcd\|abcgd | imabcd\|abcgd | imabcd\|abcgd | imabcd\|abcgd |
| be120.8.8 | 121 ; 21420 \| 121 | 5979 \| 35240 | 9.9-7 \| 9.9-7 | -2.8-7 \| -6.4-8 | 1:01 \| 4:52 |
| be120.8.9 | 121 ; 21420 \| 121 | 5788 \| 37397 | 9.9-7 \| 9.9-7 | -3.3-7 \| -8.5-8 | 1:02 \| 4:42 |
| be120.8.10 | 121 ; 21420 \| 121 | 5636 \| 35274 | 9.9-7 \| 9.9-7 | -3.5-7 \| -8.0-8 | 58 \| 4:53 |
| be250.1 | 251 ; 93375 \| 251 | 3958 \| 25038 | 9.9-7 \| 9.9-7 | 1.4-7 \| -4.9-8 | 2:10 \| 9:40 |
| be250.2 | 251 ; 93375 \| 251 | 4213 \| 29313 | 9.9-7 \| 9.9-7 | -3.7-7 \| -6.8-8 | 2:22 \| 11:36 |
| be250.3 | 251 ; 93375 \| 251 | 4230 \| 27211 | 9.9-7 \| 9.8-7 | -3.7-7 \| -4.4-8 | 2:29 \| 10:56 |
| be250.4 | 251 ; 93375 \| 251 | 4059 \| 28985 | 9.9-7 \| 9.9-7 | -3.6-7 \| -5.8-8 | 2:24 \| 11:15 |
| be250.5 | 251 ; 93375 \| 251 | 4361 \| 29277 | 9.9-7 \| 9.9-7 | -3.9-7 \| -5.2-8 | 2:35 \| 11:45 |
| bqp100-1 | 101 ; 14850 \| 101 | 7344 \| 50000 | 9.9-7 \| 1.1-6 | -9.8-8 \| -1.1-7 | 1:08 \| 5:25 |
| bqp100-2 | 101 ; 14850 \| 101 | 3799 \| 24170 | 9.9-7 \| 9.9-7 | -1.5-7 \| -6.7-8 | 30 \| 2:48 |
| bqp100-3 | 101 ; 14850 \| 101 | 3630 \| 22570 | 9.9-7 \| 9.9-7 | 8.6-8 \| -5.1-8 | 29 \| 2:35 |
| bqp100-4 | 101 ; 14850 \| 101 | 4293 \| 27893 | 9.9-7 \| 9.9-7 | -2.2-7 \| -5.9-8 | 35 \| 3:13 |
| bqp100-5 | 101 ; 14850 \| 101 | 5145 \| 34243 | 9.9-7 \| 9.9-7 | -1.0-7 \| -4.8-8 | 43 \| 3:29 |
| bqp500-1 | 501 ; 374250 \| 501 | 6385 \| 50000 | 9.9-7 \| 1.3-6 | -1.2-6 \| -1.2-7 | 23:40 \| 1:43:49 |
| bqp500-2 | 501 ; 374250 \| 501 | 6622 \| 50000 | 9.9-7 \| 1.7-6 | -1.1-6 \| -1.6-7 | 23:21 \| 1:43:49 |
| bqp500-3 | 501 ; 374250 \| 501 | 6042 \| 50000 | 9.9-7 \| 1.1-6 | -1.1-6 \| -9.1-8 | 22:10 \| 1:45:49 |
| bqp500-4 | 501 ; 374250 \| 501 | 5537 \| 50000 | 9.9-7 \| 1.2-6 | -1.1-6 \| -8.0-8 | 20:05 \| 1:46:16 |
| gka1e | 201 ; 59700 \| 201 | 5292 \| 37861 | 9.9-7 \| 9.9-7 | -2.6-7 \| -4.8-8 | 2:05 \| 10:59 |
| gka2e | 201 ; 59700 \| 201 | 4623 \| 29338 | 9.9-7 \| 9.9-7 | -6.8-7 \| -7.1-8 | 1:47 \| 8:29 |
| gka3e | 201 ; 59700 \| 201 | 6033 \| 40016 | 9.9-7 \| 9.9-7 | -3.7-7 \| -6.0-8 | 2:12 \| 11:44 |
| gka4e | 201 ; 59700 \| 201 | 7001 \| 47779 | 9.9-7 \| 9.9-7 | -5.9-7 \| -6.9-8 | 2:45 \| 14:09 |
| gka5e | 201 ; 59700 \| 201 | 6245 \| 42175 | 9.9-7 \| 9.9-7 | -5.3-7 \| -7.8-8 | 2:30 \| 12:30 |

Table 1: The performance of imABCD and ABCGD with accuracy $\varepsilon = 10^{-6}$.

We also compare our imABCD algorithm with some other BCD-type methods. The first one is a direct four-block BCD method. In this case, the block $z$ is solved by the APG-SNCG algorithm, while other blocks have analytical solutions. The second one is an enhanced version of the four-block inexact randomized ABCD method (denoted as eRABCD) that is modified from [23], where we use the proximal terms $\frac{1}{2}\|y-y^k\|^2_{\mathcal{A}\mathcal{A}^*}$ instead of $\frac{1}{2}\|y-y^k\|^2_{\lambda_{\max}(\mathcal{A}\mathcal{A}^*)}$ when updating the block $y^{k+1}$, and $\frac{1}{2}\|z-z^k\|^2_{\mathcal{B}\mathcal{B}^*+\|\mathcal{B}\|\mathcal{I}}$ when updating the block $z^{k+1}$. A similar modification has also been used in [39] when the randomized BCD algorithm is used to solve a class of positive semidefinite feasibility problems. The detailed steps of the eRABCD are given below.

---

**eRABCD:** A four-block inexact **e**nhanced **R**andomized **ABCD** algorithm for solving (19)

---

**Initialization**. Choose an initial point $W^1 = \widetilde{W}^0 \in \mathbb{W}$. Set $k = 1$ and $\alpha_0 = \frac{1}{4}$. Let $\{\varepsilon_k\}$ be a given summable sequence of error tolerance such that the error vector $\delta_z^k \in \mathbb{R}^{m_I}$ satisfies $\|\delta_z^k\| \leq \varepsilon_k$.

**Step 1.** Compute $\alpha^k = \dfrac{1}{2}\left(\sqrt{\alpha_{k-1}^4 + 4\alpha_{k-1}^2} - \alpha_{k-1}^2\right)$.

**Step 2.** Compute $\widehat{W}^{k+1} = (1-\alpha_k)\widehat{W}^k + \alpha_k \widetilde{W}^k$.

**Step 3.** Denote $\widehat{R}^k = \mathcal{A}^*\hat{y}^k + \mathcal{B}^*\hat{z}^k + \widehat{S}^k + \widehat{Z}^k + G$. Choose $i_k \in \{1,2,3,4\}$ uniformly at random and update $\widetilde{W}^{k+1}_{i_k}$ according to the following rule if the $k$-th block is selected:

$$\begin{cases} i_k = 1: & \tilde{y}^{k+1} = (\mathcal{A}\mathcal{A}^*)^{-1}((b - \mathcal{A}\widehat{R}^k)/(4\alpha_k) + \mathcal{A}\mathcal{A}^*\tilde{y}_E^k), \\ i_k = 2: & \tilde{z}^{k+1} \approx \underset{z\geq 0}{\operatorname{argmin}} \left\{ \langle \nabla_z h(\widehat{W}^{k+1}), z\rangle + \frac{4\alpha_k}{2}\|z - \tilde{z}^k\|^2_{\mathcal{B}\mathcal{B}^* + \|\mathcal{B}\|\mathcal{I}} + \langle z, \delta_z^k\rangle \right\}, \\ i_k = 3: & \widetilde{Z}^{k+1} = \Pi_{\geq 0}(\widetilde{Z}^k - \widehat{R}^k/(4\alpha_k)), \\ i_k = 4: & \widetilde{S}^{k+1} = \Pi_{\mathbb{S}^n_+}(\widetilde{S}^k - \widehat{R}^k/(4\alpha_k)). \end{cases}$$

Set $\widetilde{W}_i^{k+1} = \widetilde{W}_i^k$ for all $i \neq i_k$, $k = 1, 2, 3, 4$.

**Step 4.** Set $W_i^{k+1} = \begin{cases} \widehat{W}_i^k + 4\alpha_k(\widetilde{W}_i^{k+1} - \widetilde{W}_i^k), & i = i_k, \\ \widehat{W}_i^k, & i \neq i_k, \end{cases} \quad i = 1, 2, 3, 4.$

---

In order to know whether our proposed APG-SNCG method could universally improve the efficiency for different outer loops, we also test two variants of the BCD and eRABCD, where the block $z$ is updated by the proximal gradient step. They are named as mBCD and eRABCD2. The numerical performance of two selected test examples is shown in Table 2. One can see that the mBCD and eRABCD2 perform much worse than their inexact counterparts. These numerical results may indicate that if one of the blocks is computationally intensive (such as the block $S$ in (19) that requires the eigenvalue decomposition for each update), a small proximal term is always preferred for the other blocks in order to reduce the the number of iterations taken by the algorithm, which is also the number of updates required by the difficult block.

| | iteration | $\eta$ | time |
|---|---|---|---|
| problem | erabcd\|erabcd2\|bcd\|mbcd | erabcd\|erabcd2\|bcd\|mbcd | erabcd\|erabcd2\|bcd\|mbcd |
| bqp50-1 | 5168 \| 20172 \| 110848 \| 500000 | 9.9-6 \| 9.9-6 \| 9.9-6 \| **7.6-3** | 19 \| 40 \| 8:02 \| 18:02 |
| bqp100-1 | 7789 \| 39167 \| 203733 \| 500000 | 9.3-6 \| 9.0-6 \| 9.9-6 \| **1.1-2** | 48 \| 2:19 \| 33:24 \| 52:57 |

Table 2: The performance of eRABCD, eRABCD2, BCD and mBCD with accuracy $\varepsilon = 10^{-5}$.

Table 3 lists the numerical performance of the imABCD, BCD and eRABCD methods, with the performance profile given in Figure 2. One can see that the BCD algorithm is much less efficient than the other algorithms, as all the test examples cannot be solved to the required accuracy within 50000 iteration steps (We thus do not include its performance in the performance profile). This phenomenon has clearly demonstrated the power of the acceleration technique. Observe that the imABCD method is about 3.5 times faster than the eRABCD method.

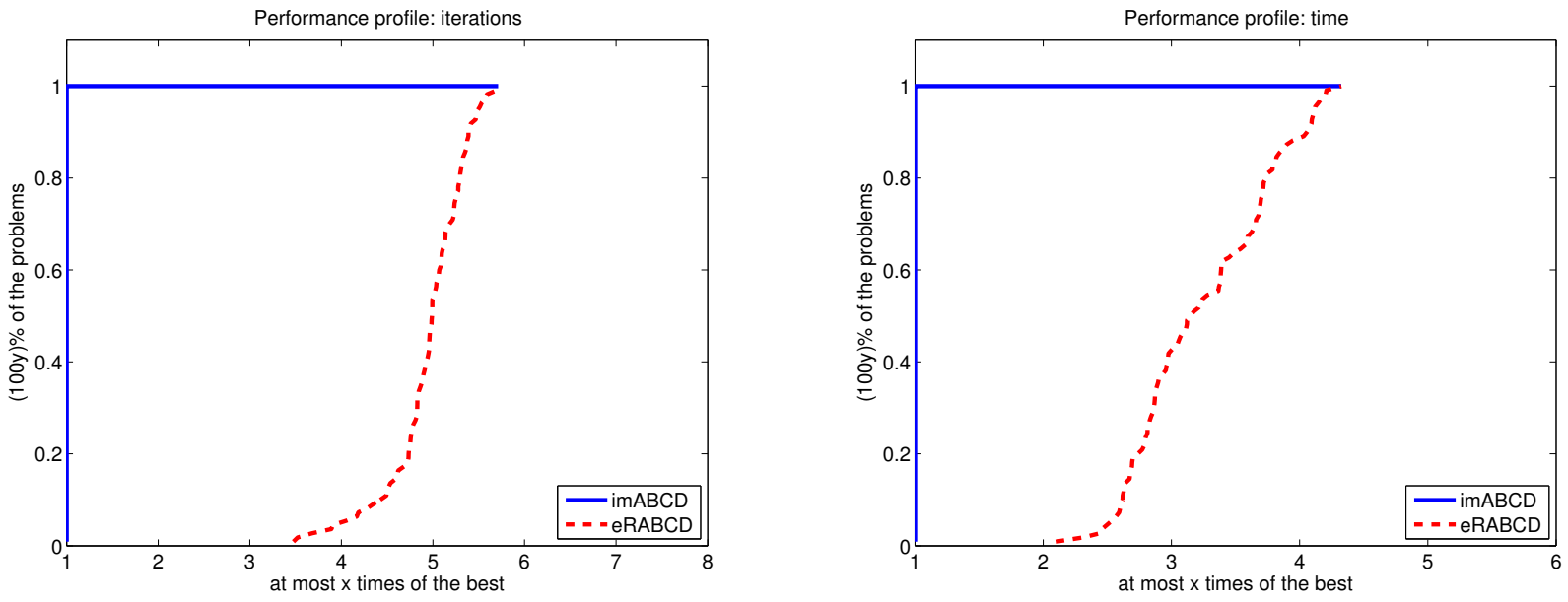

Figure 2: Performance profile of imABCD and eRABCD for with accuracy $\varepsilon = 10^{-6}$.

| | | iterations | $\eta$ | $\eta_{\text{gap}}$ | time |
|---|---|---|---|---|---|
| problem | $m_E; m_I$ \| $n_s$ | imabcd\|bcd\|erabcd | imabcd\|bcd\|erabcd | imabcd\|bcd\|erabcgd | imabcd\|bcd\|erabcd |
| be100.1 | 101 ; 14850 \| 101 | 5276 \| 50000 \| 27844 | 9.9-7 \| 3.3-5 \| 9.9-7 | -2.3-7 \| -6.1-6 \| -4.8-7 | 45 \|7:52 \|2:39 |
| be100.2 | 101 ; 14850 \| 101 | 5747 \| 50000 \| 28405 | 9.9-7 \| 3.6-5 \| 9.9-7 | -2.0-7 \| -7.6-6 \| -7.1-7 | 47 \|7:50 \|2:40 |
| be100.3 | 101 ; 14850 \| 101 | 5950 \| 50000 \| 30325 | 9.9-7 \| 3.2-5 \| 9.9-7 | -3.7-7 \| -5.1-6 \| -7.1-7 | 53 \|8:00 \|3:00 |
| be100.4 | 101 ; 14850 \| 101 | 5704 \| 50000 \| 28722 | 9.9-7 \| 3.7-5 \| 9.9-7 | -2.7-7 \| -6.7-6 \| -5.9-7 | 49 \|7:49 \|2:47 |
| be100.5 | 101 ; 14850 \| 101 | 5762 \| 50000 \| 29125 | 9.9-7 \| 3.1-5 \| 9.9-7 | -3.1-7 \| -4.5-6 \| -7.0-7 | 51 \|7:55 \|2:52 |
| be100.6 | 101 ; 14850 \| 101 | 5769 \| 50000 \| 27866 | 9.9-7 \| 3.3-5 \| 9.9-7 | -1.8-7 \| -5.9-6 \| -5.0-7 | 48 \|7:49 \|2:44 |
| be100.7 | 101 ; 14850 \| 101 | 4994 \| 50000 \| 27913 | 9.9-7 \| 3.4-5 \| 9.7-7 | -3.3-7 \| -8.9-6 \| -5.5-7 | 40 \|7:40 \|2:46 |
| be100.8 | 101 ; 14850 \| 101 | 5613 \| 50000 \| 28035 | 9.9-7 \| 3.5-5 \| 9.9-7 | -1.5-7 \| -6.7-6 \| -5.4-7 | 45 \|7:46 \|2:44 |
| be100.9 | 101 ; 14850 \| 101 | 5763 \| 50000 \| 28728 | 9.9-7 \| 3.5-5 \| 9.9-7 | -2.7-7 \| -4.0-6 \| -4.5-7 | 51 \|8:08 \|2:52 |
| be100.10 | 101 ; 14850 \| 101 | 5260 \| 50000 \| 27868 | 9.9-7 \| 3.3-5 \| 9.9-7 | -3.7-7 \| -6.7-6 \| -6.9-7 | 43 \|7:49 \|2:41 |
| be120.3.1 | 121 ; 21420 \| 121 | 4120 \| 50000 \| 22564 | 9.9-7 \| 3.6-5 \| 7.7-7 | -2.2-7 \| -1.0-5 \| -6.5-7 | 41 \|10:03 \|2:45 |
| be120.3.2 | 121 ; 21420 \| 121 | 4106 \| 50000 \| 20388 | 9.9-7 \| 3.7-5 \| 9.9-7 | -2.3-7 \| -1.1-5 \| -9.9-7 | 40 \|10:02 \|2:29 |
| be120.3.3 | 121 ; 21420 \| 121 | 3548 \| 50000 \| 18503 | 9.9-7 \| 3.5-5 \| 9.9-7 | -9.8-8 \| -1.2-5 \| -9.9-7 | 35 \|9:59 \|2:15 |
| be120.3.4 | 121 ; 21420 \| 121 | 4745 \| 50000 \| 24812 | 9.9-7 \| 3.7-5 \| 9.9-7 | -3.1-7 \| -1.2-5 \| -7.3-7 | 47 \|10:04 \|3:01 |
| be120.3.5 | 121 ; 21420 \| 121 | 5637 \| 50000 \| 27721 | 9.9-7 \| 3.6-5 \| 9.9-7 | -4.9-8 \| -8.9-6 \| -6.2-7 | 58 \|10:19 \|3:28 |
| be120.3.6 | 121 ; 21420 \| 121 | 3946 \| 50000 \| 18775 | 9.9-7 \| 3.3-5 \| 9.8-7 | -1.8-7 \| -9.0-6 \| -9.4-7 | 39 \|10:05 \|2:14 |
| be120.3.7 | 121 ; 21420 \| 121 | 4169 \| 50000 \| 22564 | 9.9-7 \| 3.5-5 \| 7.2-7 | -2.7-7 \| -1.2-5 \| -6.2-7 | 41 \|10:05 \|2:48 |
| be120.3.8 | 121 ; 21420 \| 121 | 3793 \| 50000 \| 20388 | 9.9-7 \| 3.3-5 \| 9.9-7 | -1.6-7 \| -1.2-5 \| -7.2-7 | 35 \|9:56 \|2:27 |
| be120.3.9 | 121 ; 21420 \| 121 | 4951 \| 50000 \| 22879 | 9.9-7 \| 3.7-5 \| 9.9-7 | -2.0-7 \| -1.0-5 \| -5.0-7 | 52 \|10:17 \|2:32 |
| be120.3.10 | 121 ; 21420 \| 121 | 4264 \| 50000 \| 22877 | 9.9-7 \| 3.3-5 \| 9.9-7 | -3.4-7 \| -9.8-6 \| -5.6-7 | 42 \|10:33 \|2:57 |
| be120.8.1 | 121 ; 21420 \| 121 | 5671 \| 50000 \| 28844 | 9.9-7 \| 3.6-5 \| 9.9-7 | -3.5-7 \| -9.8-6 \| -1.0-6 | 58 \|11:02 \|3:43 |
| be120.8.2 | 121 ; 21420 \| 121 | 5897 \| 50000 \| 28010 | 9.9-7 \| 3.8-5 \| 9.9-7 | -3.1-7 \| -7.1-6 \| -8.1-7 | 1:02 \|11:09 \|3:37 |
| be120.8.3 | 121 ; 21420 \| 121 | 5199 \| 50000 \| 28010 | 9.9-7 \| 3.7-5 \| 9.9-7 | -4.9-7 \| -8.5-6 \| -9.0-7 | 52 \|11:02 \|3:35 |
| be120.8.4 | 121 ; 21420 \| 121 | 6688 \| 50000 \| 33698 | 9.9-7 \| 3.6-5 \| 9.9-7 | -4.0-7 \| -6.3-6 \| -7.0-7 | 1:12 \|11:00 \|4:27 |
| be120.8.5 | 121 ; 21420 \| 121 | 5828 \| 50000 \| 28138 | 9.9-7 \| 3.6-5 \| 9.9-7 | -3.6-7 \| -5.8-6 \| -7.5-7 | 1:03 \|11:05 \|3:42 |
| be120.8.6 | 121 ; 21420 \| 121 | 4735 \| 50000 \| 25034 | 9.9-7 \| 3.6-5 \| 9.9-7 | -4.8-7 \| -9.8-6 \| -8.7-7 | 47 \|10:54 \|3:12 |
| be120.8.7 | 121 ; 21420 \| 121 | 4456 \| 50000 \| 23393 | 9.9-7 \| 3.4-5 \| 9.6-7 | -3.8-7 \| -9.3-6 \| -7.8-7 | 44 \|10:44 \|3:01 |
| be120.8.8 | 121 ; 21420 \| 121 | 5979 \| 50000 \| 28880 | 9.9-7 \| 3.7-5 \| 9.9-7 | -2.8-7 \| -8.9-6 \| -9.7-7 | 1:01 \|10:57 \|3:44 |
| be120.8.9 | 121 ; 21420 \| 121 | 5788 \| 50000 \| 28851 | 9.9-7 \| 3.8-5 \| 9.9-7 | -3.3-7 \| -6.1-6 \| -8.2-7 | 1:02 \|11:16 \|3:48 |
| be120.8.10 | 121 ; 21420 \| 121 | 5636 \| 50000 \| 28138 | 9.9-7 \| 3.6-5 \| 9.6-7 | -3.5-7 \| -7.7-6 \| -9.2-7 | 58 \|11:15 \|3:45 |
| be150.3.1 | 151 ; 33525 \| 151 | 5917 \| 50000 \| 28133 | 9.9-7 \| 3.9-5 \| 9.8-7 | -1.8-7 \| -1.2-5 \| -9.3-7 | 1:23 \|15:44 \|5:08 |
| be150.3.2 | 151 ; 33525 \| 151 | 4565 \| 50000 \| 23958 | 9.9-7 \| 3.9-5 \| 9.7-7 | -5.8-7 \| -1.4-5 \| -1.2-6 | 1:02 \|15:32 \|4:23 |
| be150.3.3 | 151 ; 33525 \| 151 | 4524 \| 50000 \| 23103 | 9.9-7 \| 4.0-5 \| 9.5-7 | -4.9-7 \| -1.3-5 \| -8.4-7 | 1:02 \|15:34 \|4:09 |
| be150.3.4 | 151 ; 33525 \| 151 | 4506 \| 50000 \| 22609 | 9.9-7 \| 3.7-5 \| 9.3-7 | -4.7-7 \| -1.3-5 \| -1.1-6 | 1:01 \|15:34 \|4:02 |
| be150.3.5 | 151 ; 33525 \| 151 | 5975 \| 50000 \| 28854 | 9.9-7 \| 3.8-5 \| 9.7-7 | -3.0-7 \| -9.8-6 \| -9.9-7 | 1:24 \|15:53 \|5:12 |
| be150.3.6 | 151 ; 33525 \| 151 | 4131 \| 50000 \| 23008 | 9.9-7 \| 3.9-5 \| 9.9-7 | -2.8-7 \| -1.7-5 \| -6.6-7 | 56 \|16:56 \|4:03 |
| be150.3.7 | 151 ; 33525 \| 151 | 4914 \| 50000 \| 23971 | 9.9-7 \| 4.0-5 \| 9.9-7 | -4.9-7 \| -1.1-5 \| -1.2-6 | 1:08 \|15:41 \|4:20 |
| be150.3.8 | 151 ; 33525 \| 151 | 4354 \| 50000 \| 22976 | 9.9-7 \| 3.6-5 \| 9.9-7 | -4.0-7 \| -1.4-5 \| -8.0-7 | 59 \|15:30 \|4:03 |
| be150.3.9 | 151 ; 33525 \| 151 | 5559 \| 50000 \| 25534 | 9.9-7 \| 3.9-5 \| 9.9-7 | -2.6-7 \| -1.1-5 \| -1.0-6 | 1:22 \|15:46 \|4:38 |
| be150.3.10 | 151 ; 33525 \| 151 | 5647 \| 50000 \| 28752 | 9.9-7 \| 3.9-5 \| 9.9-7 | -2.8-7 \| -1.2-5 \| -7.0-7 | 1:24 \|15:35 \|5:10 |
| be150.8.1 | 151 ; 33525 \| 151 | 5999 \| 50000 \| 29705 | 9.9-7 \| 3.7-5 \| 9.9-7 | -4.0-7 \| -9.0-6 \| -1.0-6 | 1:26 \|15:48 \|5:19 |
| be150.8.2 | 151 ; 33525 \| 151 | 6516 \| 50000 \| 33468 | 9.9-7 \| 3.7-5 \| 9.8-7 | -5.8-7 \| -8.3-6 \| -9.7-7 | 1:36 \|16:18 \|6:03 |
| be150.8.3 | 151 ; 33525 \| 151 | 6733 \| 50000 \| 32219 | 9.9-7 \| 4.0-5 \| 9.8-7 | -5.5-7 \| -8.3-6 \| -9.9-7 | 1:40 \|16:10 \|6:00 |
| be150.8.4 | 151 ; 33525 \| 151 | 6063 \| 50000 \| 29344 | 9.9-7 \| 3.8-5 \| 9.8-7 | -4.6-7 \| -8.0-6 \| -1.1-6 | 1:27 \|16:05 \|5:17 |
| be150.8.5 | 151 ; 33525 \| 151 | 6717 \| 50000 \| 33518 | 9.9-7 \| 4.2-5 \| 9.9-7 | -4.9-7 \| -7.6-6 \| -8.2-7 | 1:44 \|16:20 \|6:18 |
| be150.8.6 | 151 ; 33525 \| 151 | 5953 \| 50000 \| 29344 | 9.9-7 \| 3.6-5 \| 9.8-7 | -3.4-7 \| -9.7-6 \| -1.2-6 | 1:25 \|15:50 \|5:13 |
| be150.8.7 | 151 ; 33525 \| 151 | 5965 \| 50000 \| 31697 | 9.9-7 \| 3.7-5 \| 9.9-7 | -4.7-7 \| -1.0-5 \| -1.2-6 | 1:25 \|15:44 \|5:47 |
| be150.8.8 | 151 ; 33525 \| 151 | 6396 \| 50000 \| 33130 | 9.9-7 \| 3.9-5 \| 8.0-7 | -5.3-7 \| -7.9-6 \| -9.0-7 | 1:38 \|16:21 \|6:09 |
| be150.8.9 | 151 ; 33525 \| 151 | 5892 \| 50000 \| 28854 | 9.9-7 \| 4.0-5 \| 9.5-7 | -4.0-7 \| -7.4-6 \| -8.1-7 | 1:25 \|16:15 \|5:19 |

| | | | iterations | | | $\eta$ | | | $\eta_{\text{gap}}$ | | | time | | |
|---|---|---|---|---|---|---|---|---|---|---|---|---|---|---|
| problem | $m_E; m_I$ | $n_s$ | imabcd | bcd | erabcd | imabcd | bcd | erabcd | imabcd | bcd | erabcgd | imabcd | bcd | erabcd |
| be150.8.10 | 151 ; 33525 | 151 | 6262 | 50000 | 29344 | 9.9-7 | 4.0-5 | 9.7-7 | -3.7-7 | -9.4-6 | -1.1-6 | 1:30 | 16:08 | 4:24 |
| be200.3.1 | 201 ; 59700 | 201 | 6020 | 50000 | 28718 | 9.9-7 | 4.1-5 | 9.7-7 | -4.1-7 | -1.5-5 | -1.0-6 | 2:27 | 23:30 | 7:18 |
| be200.3.2 | 201 ; 59700 | 201 | 5440 | 50000 | 28026 | 9.9-7 | 4.3-5 | 9.8-7 | -6.4-7 | -1.6-5 | -1.2-6 | 2:13 | 25:08 | 6:23 |
| be200.3.3 | 201 ; 59700 | 201 | 6150 | 50000 | 30498 | 9.9-7 | 4.3-5 | 9.8-7 | -4.6-7 | -1.3-5 | -1.4-6 | 2:31 | 22:09 | 6:46 |
| be200.3.4 | 201 ; 59700 | 201 | 6079 | 50000 | 28715 | 9.9-7 | 4.4-5 | 9.8-7 | -4.1-7 | -1.3-5 | -1.3-6 | 2:30 | 22:10 | 6:16 |
| be200.3.5 | 201 ; 59700 | 201 | 6536 | 50000 | 33573 | 9.9-7 | 4.1-5 | 9.9-7 | -5.8-7 | -1.2-5 | -1.3-6 | 2:42 | 22:04 | 7:43 |
| be200.3.6 | 201 ; 59700 | 201 | 5193 | 50000 | 28028 | 9.9-7 | 4.3-5 | 9.9-7 | -5.0-7 | -1.6-5 | -1.1-6 | 2:08 | 24:24 | 6:07 |
| be200.3.7 | 201 ; 59700 | 201 | 5750 | 50000 | 28715 | 9.9-7 | 4.2-5 | 9.8-7 | -1.7-7 | -1.6-5 | -1.1-6 | 2:22 | 23:04 | 6:19 |
| be200.3.8 | 201 ; 59700 | 201 | 6087 | 50000 | 28820 | 9.9-7 | 4.1-5 | 9.9-7 | -5.2-7 | -1.4-5 | -1.5-6 | 2:30 | 22:00 | 6:35 |
| be200.3.9 | 201 ; 59700 | 201 | 6088 | 50000 | 28820 | 9.9-7 | 4.2-5 | 9.9-7 | -3.9-7 | -1.4-5 | -1.4-6 | 2:32 | 22:16 | 6:29 |
| be200.3.10 | 201 ; 59700 | 201 | 5464 | 50000 | 28026 | 9.9-7 | 4.0-5 | 9.6-7 | -6.8-7 | -1.6-5 | -1.3-6 | 2:14 | 22:14 | 6:16 |
| be200.8.1 | 201 ; 59700 | 201 | 6882 | 50000 | 34136 | 9.9-7 | 4.0-5 | 9.9-7 | -7.2-7 | -7.9-6 | -1.5-6 | 2:53 | 23:08 | 8:00 |
| be200.8.2 | 201 ; 59700 | 201 | 6345 | 50000 | 33798 | 9.9-7 | 4.0-5 | 9.8-7 | -7.0-7 | -1.0-5 | -1.4-6 | 2:37 | 22:44 | 7:48 |
| be200.8.3 | 201 ; 59700 | 201 | 7829 | 50000 | 39097 | 9.9-7 | 4.1-5 | 9.3-7 | -6.4-7 | -8.3-6 | -1.0-6 | 3:17 | 23:29 | 9:09 |
| be200.8.4 | 201 ; 59700 | 201 | 7475 | 50000 | 34483 | 9.9-7 | 4.2-5 | 9.8-7 | -6.3-7 | -6.7-6 | -1.2-6 | 3:21 | 24:03 | 8:15 |
| be200.8.5 | 201 ; 59700 | 201 | 7600 | 50000 | 39771 | 9.8-7 | 4.1-5 | 9.9-7 | -7.1-7 | -6.8-6 | -9.9-7 | 3:24 | 23:58 | 9:40 |
| be200.8.6 | 201 ; 59700 | 201 | 7164 | 50000 | 34143 | 9.9-7 | 3.9-5 | 9.9-7 | -7.6-7 | -7.6-6 | -1.6-6 | 3:00 | 23:11 | 7:47 |
| be200.8.7 | 201 ; 59700 | 201 | 5677 | 50000 | 29693 | 9.9-7 | 3.9-5 | 9.7-7 | -7.4-7 | -1.5-5 | -1.9-6 | 2:20 | 22:18 | 5:42 |
| be200.8.8 | 201 ; 59700 | 201 | 5620 | 50000 | 29427 | 9.9-7 | 3.9-5 | 9.9-7 | -7.7-7 | -1.1-5 | -1.7-6 | 2:19 | 21:33 | 6:32 |
| be200.8.9 | 201 ; 59700 | 201 | 7100 | 50000 | 34143 | 9.9-7 | 4.0-5 | 9.9-7 | -7.1-7 | -7.4-6 | -1.5-6 | 3:01 | 22:11 | 7:38 |
| be200.8.10 | 201 ; 59700 | 201 | 6981 | 50000 | 34136 | 9.9-7 | 3.9-5 | 9.7-7 | -7.5-7 | -8.3-6 | -1.5-6 | 2:55 | 21:58 | 7:35 |
| be250.1 | 251 ; 93375 | 251 | 3958 | 50000 | 19178 | 9.9-7 | 5.2-5 | 9.8-7 | 1.4-7 | -3.0-5 | -1.9-6 | 2:10 | 35:43 | 6:15 |
| be250.2 | 251 ; 93375 | 251 | 4213 | 50000 | 23028 | 9.9-7 | 5.1-5 | 9.9-7 | -3.7-7 | -2.8-5 | -1.1-6 | 2:22 | 37:47 | 7:25 |
| be250.3 | 251 ; 93375 | 251 | 4230 | 50000 | 23195 | 9.9-7 | 5.1-5 | 9.8-7 | -3.7-7 | -3.0-5 | -1.2-6 | 2:29 | 37:37 | 7:52 |
| be250.4 | 251 ; 93375 | 251 | 4059 | 50000 | 23109 | 9.9-7 | 4.7-5 | 9.9-7 | -3.6-7 | -2.4-5 | -1.1-6 | 2:24 | 32:53 | 7:45 |
| be250.5 | 251 ; 93375 | 251 | 4361 | 50000 | 23141 | 9.9-7 | 5.2-5 | 9.9-7 | -3.9-7 | -2.7-5 | -1.1-6 | 2:35 | 38:43 | 7:39 |
| be250.6 | 251 ; 93375 | 251 | 4131 | 50000 | 22773 | 9.9-7 | 5.2-5 | 7.7-7 | -2.4-7 | -2.9-5 | -8.7-7 | 2:25 | 37:40 | 7:46 |
| be250.7 | 251 ; 93375 | 251 | 6019 | 50000 | 28518 | 9.9-7 | 4.9-5 | 9.9-7 | -1.8-7 | -2.1-5 | -1.1-6 | 3:37 | 33:21 | 9:46 |
| be250.8 | 251 ; 93375 | 251 | 3398 | 50000 | 18307 | 9.9-7 | 5.5-5 | 8.6-7 | -5.0-7 | -3.3-5 | -1.4-6 | 2:00 | 38:14 | 5:56 |
| be250.9 | 251 ; 93375 | 251 | 5299 | 50000 | 28464 | 9.9-7 | 4.9-5 | 9.9-7 | -5.8-7 | -2.4-5 | -1.0-6 | 3:09 | 33:24 | 9:38 |
| be250.10 | 251 ; 93375 | 251 | 4601 | 50000 | 24611 | 9.9-7 | 5.0-5 | 9.9-7 | -3.7-7 | -2.6-5 | -1.1-6 | 2:39 | 37:58 | 8:09 |
| bqp50-1 | 51 ; 3675 | 51 | 3219 | 50000 | 13417 | 9.9-7 | 2.3-5 | 9.8-7 | -2.1-8 | -2.8-6 | -5.3-8 | 16 | 4:03 | 49 |
| bqp50-2 | 51 ; 3675 | 51 | 4717 | 50000 | 19129 | 9.9-7 | 3.2-5 | 9.9-7 | -3.7-8 | -3.5-6 | -1.3-7 | 21 | 3:52 | 1:00 |
| bqp50-3 | 51 ; 3675 | 51 | 3252 | 50000 | 14389 | 9.9-7 | 2.9-5 | 9.9-7 | -6.4-8 | -3.8-6 | -1.7-7 | 14 | 3:49 | 44 |
| bqp50-4 | 51 ; 3675 | 51 | 3817 | 50000 | 14850 | 9.9-7 | 1.9-5 | 9.9-7 | -3.6-8 | -1.3-6 | -5.7-8 | 18 | 3:54 | 51 |
| bqp50-5 | 51 ; 3675 | 51 | 1887 | 50000 | 6945 | 9.8-7 | 1.2-5 | 9.9-7 | 3.4-8 | -1.1-6 | -7.6-8 | 09 | 4:00 | 24 |
| bqp50-6 | 51 ; 3675 | 51 | 3234 | 50000 | 14601 | 9.9-7 | 3.2-5 | 9.9-7 | -6.0-8 | -4.8-6 | -1.4-7 | 14 | 3:50 | 46 |
| bqp50-7 | 51 ; 3675 | 51 | 2209 | 50000 | 11197 | 9.9-7 | 1.9-5 | 9.9-7 | -1.6-8 | -1.9-6 | -5.9-8 | 10 | 3:51 | 36 |
| bqp50-8 | 51 ; 3675 | 51 | 4293 | 50000 | 18393 | 9.9-7 | 3.1-5 | 9.9-7 | -8.2-8 | -5.5-6 | -1.1-7 | 19 | 3:52 | 56 |
| bqp50-9 | 51 ; 3675 | 51 | 4751 | 50000 | 18573 | 9.9-7 | 3.1-5 | 9.9-7 | -1.0-8 | -4.0-6 | -8.5-8 | 22 | 3:55 | 1:04 |
| bqp50-10 | 51 ; 3675 | 51 | 4145 | 50000 | 14601 | 9.9-7 | 3.2-5 | 9.9-7 | -3.6-8 | -3.1-6 | -8.0-8 | 20 | 3:57 | 52 |
| bqp100-1 | 101 ; 14850 | 101 | 7344 | 50000 | 25527 | 9.9-7 | 4.8-5 | 9.9-7 | -9.8-8 | -6.1-6 | -2.7-7 | 1:08 | 7:54 | 2:22 |
| bqp100-2 | 101 ; 14850 | 101 | 3799 | 50000 | 18857 | 9.9-7 | 4.3-5 | 9.8-7 | -1.5-7 | -1.4-5 | -6.4-7 | 30 | 7:22 | 1:39 |
| bqp100-3 | 101 ; 14850 | 101 | 3630 | 50000 | 18066 | 9.9-7 | 4.4-5 | 6.7-7 | 8.6-8 | -1.1-5 | -2.5-7 | 29 | 7:22 | 1:32 |
| bqp100-4 | 101 ; 14850 | 101 | 4293 | 50000 | 22744 | 9.9-7 | 4.4-5 | 9.6-7 | -2.2-7 | -1.0-5 | -3.6-7 | 35 | 7:26 | 1:58 |
| bqp100-5 | 101 ; 14850 | 101 | 5145 | 50000 | 23136 | 9.9-7 | 4.5-5 | 9.9-7 | -1.0-7 | -1.0-5 | -3.7-7 | 43 | 7:25 | 1:58 |
| bqp100-6 | 101 ; 14850 | 101 | 4170 | 50000 | 20467 | 9.9-7 | 4.6-5 | 9.9-7 | -2.3-7 | -9.6-6 | -6.3-7 | 34 | 7:34 | 1:46 |
| bqp100-7 | 101 ; 14850 | 101 | 5697 | 50000 | 27700 | 9.9-7 | 4.7-5 | 9.9-7 | -1.4-7 | -9.3-6 | -3.7-7 | 48 | 7:28 | 2:27 |
| bqp100-8 | 101 ; 14850 | 101 | 5491 | 50000 | 22986 | 9.9-7 | 4.4-5 | 9.9-7 | -6.4-8 | -8.3-6 | -4.7-7 | 46 | 7:29 | 2:00 |

| | | iterations | $\eta$ | $\eta_{\text{gap}}$ | time |
|---|---|---|---|---|---|
| problem | $m_E; m_I \mid n_s$ | imabcd\|bcd\|erabcd | imabcd\|bcd\|erabcd | imabcd\|bcd\|erabcgd | imabcd\|bcd\|erabcd |
| bqp100-9 | 101 ; 14850 \| 101 | 5481 \| 50000 \| 24863 | 9.9-7 \| 4.1-5 \| 9.9-7 | -1.3-7 \| -6.7-6 \| -4.5-7 | 47 \|7:28 \|2:13 |
| bqp100-10 | 101 ; 14850 \| 101 | 5715 \| 50000 \| 24868 | 9.9-7 \| 4.4-5 \| 9.9-7 | -6.9-8 \| -1.0-5 \| -5.7-7 | 45 \|7:25 \|2:10 |
| bqp250-1 | 251 ; 93375 \| 251 | 7230 \| 50000 \| 34885 | 9.9-7 \| 4.7-5 \| 9.9-7 | -7.1-7 \| -1.4-5 \| -1.6-6 | 4:06 \|32:10 \|10:38 |
| bqp250-2 | 251 ; 93375 \| 251 | 6848 \| 50000 \| 34885 | 9.9-7 \| 4.8-5 \| 9.9-7 | -7.7-7 \| -1.5-5 \| -1.7-6 | 3:51 \|31:20 \|10:17 |
| bqp250-3 | 251 ; 93375 \| 251 | 5352 \| 50000 \| 28522 | 9.9-7 \| 4.8-5 \| 9.7-7 | -8.0-7 \| -2.1-5 \| -1.4-6 | 2:56 \|35:00 \|8:13 |
| bqp250-4 | 251 ; 93375 \| 251 | 6624 \| 50000 \| 33999 | 9.9-7 \| 5.0-5 \| 9.9-7 | -6.2-7 \| -1.3-5 \| -1.4-6 | 3:46 \|31:52 \|10:08 |
| bqp250-5 | 251 ; 93375 \| 251 | 6933 \| 50000 \| 34834 | 9.9-7 \| 4.7-5 \| 9.9-7 | -7.5-7 \| -1.5-5 \| -1.7-6 | 3:55 \|31:36 \|10:15 |
| bqp250-6 | 251 ; 93375 \| 251 | 6789 \| 50000 \| 34379 | 9.9-7 \| 4.9-5 \| 9.8-7 | -6.1-7 \| -1.4-5 \| -1.5-6 | 3:49 \|33:00 \|10:16 |
| bqp250-7 | 251 ; 93375 \| 251 | 5894 \| 50000 \| 33670 | 9.9-7 \| 5.0-5 \| 9.9-7 | -6.2-7 \| -2.1-5 \| -9.6-7 | 3:17 \|34:30 \|9:47 |
| bqp250-8 | 251 ; 93375 \| 251 | 7654 \| 50000 \| 34379 | 9.9-7 \| 5.0-5 \| 9.9-7 | -6.7-7 \| -1.1-5 \| -1.4-6 | 4:32 \|32:13 \|10:27 |
| bqp250-9 | 251 ; 93375 \| 251 | 6089 \| 50000 \| 33694 | 9.9-7 \| 4.8-5 \| 9.9-7 | -8.8-7 \| -1.9-5 \| -1.3-6 | 3:23 \|32:00 \|9:43 |
| bqp250-10 | 251 ; 93375 \| 251 | 7023 \| 50000 \| 34834 | 9.9-7 \| 4.6-5 \| 9.8-7 | -6.1-7 \| -1.7-5 \| -1.6-6 | 3:55 \|34:32 \|10:20 |

Table 3: The performance of imABCD, BCD and eRABCD with accuracy $\varepsilon = 10^{-6}$.

Based on the above numerical results, we may conclude that the efficiency of the imABCD algorithm can be attributed to the double acceleration procedure: the outer acceleration of the two-block coordinate descent method, and the inner acceleration by proper incorporation of the second-order information through solving the subproblems in each iteration by Newton type methods.

# 6 Conclusions

In this paper, for the purpose of overcoming the potential degeneracy of the matrix best approximation problem (2), we have proposed a two-block imABCD method with each block solved by the Newton type methods. Extensive numerical results demonstrated the efficiency and robustness of our algorithm in solving various instances of the large scale matrix approximation problems. We believe that our algorithmic framework is a powerful tool to deal with degenerate problems and might be adapted to other convex matrix optimization problems in the future.